\documentclass[a4paper,11pt,reqno]{amsart}
\usepackage[latin1]{inputenc}
\usepackage[english]{babel}
\usepackage{amssymb,amsfonts, amsmath, color}
\usepackage[margin=2.7cm]{geometry}
\usepackage{xcolor}
\usepackage{graphicx, color, enumerate}
\usepackage[latin1]{inputenc}
\usepackage[active]{srcltx}
\usepackage{tikz}
\usepackage{pgf}
\usepackage{etex}
\usepackage{verbatim}
\usepackage{tikz-3dplot}
\usepackage{pgfkeys}
\usepackage{amsmath}
\usepackage{amsfonts}
\usepackage{amssymb}
\usepackage{amsthm}
\usepackage{algpseudocode}
\usepackage{float}
\usepackage{xcolor}
\usepackage{mathrsfs}
\usepackage{import}
\usepackage{geometry}
\usepackage{fancyhdr}
\usepackage{fp}
\usepackage[colorlinks, citecolor=blue, linkcolor=blue]{hyperref}

\newtheorem{theorem}{Theorem}[section]
\newtheorem{proposition}[theorem]{Proposition}
\newtheorem{corollary}[theorem]{Corollary}
\newtheorem{lemma}[theorem]{Lemma}
\newtheorem{remark}[theorem]{Remark}
\newtheorem{definition}[theorem]{Definition}

\newcommand{\bcl}{\begin{center}}
\newcommand{\ecl}{\end{center}}
\newcommand{\brl}{\begin{right}}
\newcommand{\erl}{\end{right}}
\newcommand{\ben}{\begin{enumerate}}
\newcommand{\een}{\end{enumerate}}
\newcommand{\overliner}{\begin{array}}
\newcommand{\earr}{\end{array}}
\newcommand{\btab}{\begin{tabular}}
\newcommand{\etab}{\end{tabular}}
\newcommand{\bdoc}{\begin{document}}
\newcommand{\edoc}{\end{document}}
\newcommand{\beqy}{\begin{eqnarray}}
\newcommand{\eeqy}{\end{eqnarray}}

\newcommand{\beqi}{\begin{eqnarray*}}
\newcommand{\eeqi}{\end{eqnarray*}}
\newcommand{\bitem}{\begin{itemize}}

\newcommand{\eitem}{\end{itemize}}
\newcommand{\nln}{\newline}
\newcommand{\newt}{\newtheorem}

\newcommand{\pa}{\partial}
\newcommand{\re}{{I\!\!R}}
\newcommand{\Rn}{\R^N}
\newcommand{\xr}{x\in\R }
\newcommand{\x}{\times}
\newcommand{\dyle}{\displaystyle}
\newcommand{\ene}{{I\!\!N}}
\newcommand{\irn}{\int\limits_{\R^N}}
\newcommand{\io}{\int\limits_{\O}}
\newcommand{\meas}{{\rm meas\,}}

\newcommand{\sign}{{\rm sign}}
\newcommand{\map}{\longrightarrow }
\newcommand{\imp}{\Longrightarrow }
\renewcommand{\div}{\nabla\cdot }
\newcommand{\sen}{{\rm sen\,}}
\newcommand{\tg}{{\rm tg\,}}
\newcommand{\arcsen}{{\rm arcsen\,}}
\newcommand{\arctg}{{\rm arctg\,}}
\newcommand{\supp}{{\textsl supp\ }}
\newcommand{\ity}{\int_{-\iy}^{+\iy}}
\newcommand{\limit}{\lim\limits}
\newcommand{\limi}{\limit_{n\to\infty}}
\newcommand{\sumi}{\sum\limits_{n=1}^{\infty}}
\newcommand{\ulu}{\underline u}
\newcommand{\ulw}{\underline w}
\newcommand{\ulz}{\underline z}
\newcommand{\ulv}{\underline v}
\newcommand{\uls}{\underline s}
\newcommand{\olu}{\overline u}
\newcommand{\olv}{\overline v}
\newcommand{\ols}{\overline s}
\newcommand{\ob}{\overline\b}
\newcommand{\ovar}{\overline\var}
\newcommand{\wv}{\widetilde v}
\newcommand{\wu}{\widetilde u}
\newcommand{\ws}{\widetilde s}
\renewcommand{\a }{\alpha }
\renewcommand{\b }{\beta }
\newcommand{\g }{\gamma}
\newcommand{\G }{\Gamma }
\renewcommand{\d }{\delta }

\newcommand{\D }{\Delta }
\newcommand{\e }{\varepsilon }
\newcommand{\z }{\zeta }
\renewcommand{\l }{\lambda }
\renewcommand{\L }{\Lambda }
\newcommand{\m }{\mu }
\newcommand{\n }{\nabla }
\newcommand{\s }{\sigma }
\newcommand{\Sig }{\Sigma }
\renewcommand{\t }{\tau }
\newcommand{\var }{\varphi }
\renewcommand{\o }{\omega }
\renewcommand{\O }{\Omega }
\newcommand{\R}{{\mathbb{R}}}
\newcommand{\bC}{{\bf C}}
\newcommand{\bZ}{{\bf Z}}
\newcommand{\bN}{{\bf N}}
\newcommand{\bQ}{{\bf Q}}
\newcommand{\bK}{{\bf K}}
\newcommand{\bI}{{\bf I}}
\newcommand{\bv}{{\bf v}}
\newcommand{\bV}{{\bf V}}
\DeclareMathOperator{\suppo}{supp} \DeclareMathOperator{\di}{div}

\newenvironment{Proof}{\Rmovelastskip\vskip12pt
plus 1pt \noindent\em\rm}{\hfill {\qed \hskip .2cm}}

\begin{document}

\title[Heat equation with density on Riemannian manifolds]{Rigidity for the heat equation with density on Riemannian manifolds through a conformal change}

\author{Alexander Grigor'yan}

\address{\hbox{\parbox{5.7in}{\medskip \noindent{Alexander Grigor'yan, \\Facult\"at f\"ur Mathematik, \\Universit\"at Bielefeld, 33501, Bielefeld, Germany \\ [3pt] \emph{E-mail address: }{\tt grigor@math.uni-bielefeld.de}}}}}

\author{Giulia Meglioli}

\address{\hbox{\parbox{5.7in}{\medskip \noindent{Giulia Meglioli, \\Facult\"at f\"ur Mathematik, \\Universit\"at Bielefeld, 33501, Bielefeld, Germany \\ [3pt] \emph{E-mail address: }{\tt gmeglioli@math.uni-bielefeld.de}}}}}

\author{Alberto Roncoroni}

\address{\hbox{\parbox{5.7in}{\medskip \noindent{Alberto Roncoroni, \\Dipartimento di Matematica, \\Politecnico di Milano, \\Piazza Leonardo da Vinci 32, 20133, Milano, Italy \\ [3pt] \emph{E-mail address: }{\tt alberto.roncoroni@polimi.it}}}}}



\maketitle

\tableofcontents

\begin{abstract}
We investigate uniqueness of solution to the heat equation with a density $\rho$ on complete, non-compact weighted Riemannian manifolds of infinite volume. Our main goal is to identify sufficient conditions under which the solution $u$ vanishes identically, assuming that $u$ belongs to a certain weighted Lebesgue space with exponential or polynomial weight, $L^p_{\phi}$.
We distinguish between the cases $p > 1$ and $p = 1$ which required stronger assumptions on the manifold and the density function $\rho$. We develop a unified method based on a conformal transformation of the metric, which allows us to reduce the problem to a standard heat equation on a suitably weighted manifold. In addition, we construct explicit counterexamples on model manifolds which demonstrate optimality of our assumptions on the density $\rho$. 
\end{abstract}

\
\section{Introduction}\setcounter{equation}{0}
\

In this paper we study uniqueness phenomena for the heat equation with density posed on weighted Riemannian manifolds. More precisely, we consider a complete, non-compact Riemannian manifold $(M,g,\mu)$ of dimension $N \geq 2$ with infinite volume and we investigate uniqueness of classical solutions to the Cauchy problem
\begin{equation}\label{problema0}
\begin{cases}
\rho\,\partial_t u=\Delta u +f \quad  &\text{in}\,\,\,M\times(0,T],\\
u=u_0\quad &\text{on}\,\,\,M\times\{0\}\,.
\end{cases}
\end{equation}
Here $\Delta$ is the Laplace-Beltrami operator on $M$, $f \in C(M \times (0,T])$, $u_0 \in L^p_{\text{loc}}(M)$ for some $p \geq 1$ and \begin{equation}\label{compl}
\begin{aligned}
&\quad\quad\rho\in C^\infty(M),\,\,\,\rho>0\,\,\text{in}\,\,\,\,M\,\,\,\text{with}\\
&\int_1^\infty\sqrt{\rho(\gamma(s))}\, ds=+\infty\, ,\quad \text{ for every divergent ray $\gamma$}\, .
\end{aligned}
\end{equation}
We recall that a minimizing geodesic $\gamma:[0,+\infty)\rightarrow M$ is called a \emph{ray} and it is called \emph{divergent} if eventually it leaves every compact subset of $M$.
Clearly, uniqueness of solutions to problem \eqref{problema0} follows if one shows that $u\equiv0$ is the only solution to the following problem
\begin{equation}\label{problema}
\begin{cases}
\rho\,\partial_t u=\Delta u  \quad  &\text{in}\,\,\,M\times(0,T],\\
u=0\quad &\text{on}\,\,\,M\times\{0\}\,.
\end{cases}
\end{equation}

This framework naturally connects to the case of studying uniqueness for the heat equation of a \textit{weighted} Laplacian posed on a weighted manifolds $(M,\tilde g,\tilde\mu)$. More precisely, we can define the new volume measure $d\tilde\mu = \rho\, d\mu$ with $d\mu$ being the Riemannian volume form, and the new metric $\tilde{g} = \rho\, g$ as a conformal change of the underlying Riemannian metric $g$. This transformation induces a corresponding change in the geometric structure of the manifold, affecting both the volume form and the differential operators. In particular, the Laplace-Beltrami operator transforms accordingly, giving rise to a weighted Laplacian that reflects the influence of the density $\rho$ on the diffusion process. Under this conformal transformation, the original parabolic problem \eqref{problema} can be rewritten in terms of the conformal metric $\tilde{g}$ and the associated weighted measure $d\tilde{\mu}$ as follows:
\begin{equation}\label{introproblem}
\begin{cases}
\partial_t u = \tilde{\Delta} u & \text{in } M \times (0,T], \\
u = 0 & \text{on } M \times \{0\},
\end{cases}
\end{equation}
where $\tilde{\Delta}$ denotes the weighted Laplacian with respect to $(M, \tilde{g}, d\tilde{\mu})$ and satisfies the identity $$\tilde{\Delta} = \frac{1}{\rho} \Delta.$$ This reformulation plays a crucial role in our analysis, as it allows us to treat the problem as a classical heat equation on a weighted, conformally complete Riemannian manifold.

\subsection{Summary of known results}\label{intro}

The problem of uniqueness for parabolic equations on manifolds has been widely studied, particularly in connection with geometric analysis, stochastic completeness, and the behavior of the heat kernel. 
We begin by discussing the Euclidean setting. Let $T > 0$ and let the operator $L$ be defined by
\[
L u := \sum_{i,j=1}^N \frac{\partial^2}{\partial x_i \partial x_j} \left[a_{ij}(x,t)\, u\right] - \sum_{i=1}^N \frac{\partial}{\partial x_i} \left[b_i(x,t)\, u\right] + c(x,t)\, u,
\]
where the coefficients $a_{ij}$, $b_i$, and $c$ are locally bounded, together with all their derivatives, and the matrix $A = (a_{ij})$ is assumed to be positive semi-definite in $\mathbb{R}^N \times (0,T)$. In \cite{AB}, the author investigate uniqueness of solutions to the Cauchy problem
\begin{equation}\label{eq18}
\begin{cases}
\partial_t u = L u & \text{in } \mathbb{R}^N \times (0,T), \\
u = 0 & \text{in } \mathbb{R}^N \times \{0\}.
\end{cases}
\end{equation}
They further assume that there exist constants $\lambda \geq 0$ and $K_1, K_2, K_3 > 0$ such that, for almost every $(x,t) \in \mathbb{R}^N \times (0,T)$,
\begin{equation*}
|a_{ij}(x,t)| \leq K_1 \left(1 + |x|^2\right)^{\frac{2 - \lambda}{2}}, \quad
|b_i(x,t)| \leq K_2 \left(1 + |x|^2\right)^{\frac{1}{2}}, \quad
|c(x,t)| \leq K_3 \left(1 + |x|^2\right)^{\frac{\lambda}{2}}.
\end{equation*}
Then, in \cite[Theorem 1]{AB}, it is shown that if $u$ is a solution to problem \eqref{eq18} and belongs to $L^1_\phi(\mathbb{R}^N \times (0,T))$ with
\begin{equation}\label{eq19}
\phi(x) = (|x|^2 + 1)^{-\alpha_0} \quad \text{for } \lambda = 0,
\end{equation}
or
\begin{equation}\label{eq15}
\phi(x) = \exp\left\{ -\alpha_0 (|x|^2 + 1)^{\frac{2 - \lambda}{2}} \right\} \quad \text{for } \lambda > 0,
\end{equation}
for some $\alpha_0 > 0$, then the only solution is the trivial one
\(
u \equiv 0\) in \(\mathbb{R}^N \times (0,T).\)
Here, given a positive, continuous function $\phi \in C(\mathbb{R}^N \times (0,T))$, one has defined the weighted Lebesgue space 
\[
L^1_\phi(\mathbb{R}^N \times (0,T)) := \left\{ u : \mathbb{R}^N \times (0,T) \to \mathbb{R} \,\middle|\, \int_0^T \int_{\mathbb{R}^N} |u(x,t)| \, \phi(x,t) \, dx \, dt < \infty \right\}.
\]
The result in \cite{AB} has attracted significant attention, and similar uniqueness results have been obtained in \cite{EKP, F, IKO, KPT, PoT, Ti}, where further discussions and refinements can be found. 
Furthermore, analogous problems have been studied in the context of degenerate elliptic and parabolic equations both in $\mathbb{R}^N$ and in bounded domains; see for instance \cite{PPT2, PoT, PTes} and the references therein.

Let us now move to the Riemannian setting. Consider the classical heat equation
\begin{equation}\label{heat}
\begin{cases}
\partial_t u = \Delta u & \text{in } M \times (0,T), \\
u = u_0 & \text{on } M \times \{0\},
\end{cases}
\end{equation}
posed on a geodesically complete, non-compact Riemannian manifold $M$, i.e. \eqref{problema0} with $\rho\equiv1$, $f\equiv0$. The question of uniqueness of solutions to \eqref{heat} has been extensively studied. It is well known that uniqueness for this problem is equivalent to the stochastic completeness of the manifold $M$; see for instance \cite[Theorem 6.2]{Grig}. Given the heat kernel $p$ and the Riemannian volume element $d\mu$, we say that $M$ is stochastically complete if
$$
\int_M p(x,y,t)\,d\mu(y)=1\quad \text{for all}\,\,\,x\in M,\,\,t>0.
$$
We refer the reader, e.g., to \cite{Gaf, Grigbook, Grig, Hsu, Pinch, Pnuovo, Yau} for some results regarding the uniqueness of solutions to linear parabolic Cauchy problems related to \eqref{problema}. In particular, \cite[Theorem 9.2]{Grig} the author proves that if $u$ is a classical solution of \eqref{heat} with vanishing initial data $u_0 \equiv 0$, then
\[
u \equiv 0 \quad \text{in } M \times (0,T),
\]
provided that $u$ belongs to the weighted space
\[
L^2_\phi(M \times (0,T)) := \left\{ u : M \times (0,T) \to \mathbb{R} \,\middle|\, \int_0^T \int_M u(x,t)^2 \, \phi(x) \, d\mu(x) \, dt < \infty \right\},
\]
where $\phi(x) = e^{-f(r(x))}$, $r(x)$ is the distance from a fixed reference point in $M$, and $f$ is a continuous, positive, increasing function on $(0,\infty)$ such that, for some $R_0 > 0$
\[
\int_{R_0}^\infty \frac{r}{f(r)} \, dr = \infty.
\]
This result was later extended in \cite{Punzo2015}, where it is shown that the same uniqueness property holds for solutions $u \in L^p_\phi(M \times (0,T))$ with $1 < p \leq 2$, using the same class of weight functions $\phi$. This generalization highlights the robustness of the method and its applicability to broader classes of nonlinear or degenerate problems on manifolds. Moreover, let us mention the contribution in \cite{IM}. In this paper, the authors investigate uniqueness properties for nonnegative solutions of the Cauchy problem \eqref{problema} on noncompact Riemannian manifolds and Euclidean domains for even more general linear operators $L$ which may include some lower order terms. Let us in particular mention that they establish a sharp uniqueness result in the weighted $L^2$ space in the case of problem \eqref{problema} with $\rho$ satisfying for some $C_1,C_2>0$
$$
C_1r^{-\alpha}\left(\log r+1\right)^{-\beta}\le \rho(x)\le C_2 r^{-\alpha}\left(|\log r|+1\right)^{-\beta}
$$
where $\alpha<2$ or $\alpha=2$ and $\beta\ge1$.

More recently, similar uniqueness questions have been addressed in the discrete setting, particularly on infinite weighted graphs. In this context, problem \eqref{problema} has been studied for various choices of weights, with particular attention to the case $\rho \equiv 1$. A notable result is presented in \cite{Huang}, where uniqueness of solutions is established under a specific integrability condition. In addition, the author constructs an explicit counterexample on the integer lattice $\mathbb{Z}$ showing that his condition is somehow sharp. Another contribution in this direction is due to \cite{HKS}, where a different sufficient condition for uniqueness is provided, again for the case $\rho \equiv 1$. Furthermore, in a particular class of graphs, the authors also prove that the uniqueness of bounded solutions to problem \eqref{problema} is equivalent to the stochastic completeness of the graph. We also refer the reader to \cite{BP1,Meg,MP3} for other uniqueness investigations to problem \eqref{problema} when posed on a discrete graph.

\

We conclude the list of known result by mentioning that similar uniqueness questions have been addressed also in the case of stationary (elliptic) linear equation, see e.g. \cite{Li,LS,MR,MP2,MP3,Naber,PW,PRS2,P2,Wu,Yau2} and references therein.

\subsection{Outline of our results} 

Our main contribution is to establish uniqueness results for classical solutions to the Cauchy problem in the weighted Lebesgue spaces $L^p_{e^{-\psi}}(M \times (0,T))$, where the exponential weight $e^{-\psi}$ captures the growth behavior of the solutions at infinity. We develop a systematic framework that treats separately the cases $p > 1$ and $p = 1$, which exhibit fundamentally different analytical features and require distinct techniques.
In the case $p > 1$, we first prove a general uniqueness theorem, see Theorem \ref{teogen}, under minimal geometric assumptions. Specifically, assuming that the density $\rho$ satisfies the completeness condition \eqref{compl} and is bounded above and below by positive radial functions $\rho_1(r)$ and $\rho_2(r)$, see \eqref{eq:rhobound-general}. We show that uniqueness holds provided the weight function $\psi$ satisfies an integral condition involving both $\rho_1$ and $\rho_2$, see \eqref{integral-phi-gen}. We then specialize Theorem \ref{teogen} to two notable settings. In Theorem \ref{main1}, we consider the case where the density has polynomial decay of the form $\rho(x) \sim (1 + r^2(x))^{-\frac\theta2}$ with $0 \leq \theta < 2$, see \eqref{eq:rhobound-theta}. Here, the integral condition on $\psi$ reduces to the simpler form
\[
\int_{R_0}^\infty \frac{r^{1-\theta}}{\psi(r)} \, dr = +\infty,
\]
which allows for a direct link between the decay of $\rho$ and the admissible growth of the solution. A similar but more delicate result is obtained in Theorem \ref{main2} for the critical case $\theta = 2$, where the density decays as $\rho(x) \sim (1 + r^2(x))^{-1}$, see \eqref{eq:rhobound}. In this case, the uniqueness condition becomes
\[
\int_{R_0}^\infty \frac{\log r}{r\, \psi(r)} \, dr = +\infty,
\]
highlighting the borderline behavior of solutions near the threshold of non-uniqueness.
For the critical case $p = 1$, additional geometric structure is required. We assume that the manifold $(M,g)$ possesses a pole and that the Laplacian of the distance function from a reference point in $M$ possesses some proper lower bound. Furthermore, we restrict the analysis to radial functions $\rho$. Under these extra assumptions on the underlying manifold, we prove analogues of the previous theorems: Theorem \ref{teogen-p1} extends Theorem \ref{teogen} to the $L^1$, while Theorems \ref{main1-p1} and \ref{main2-p1} mirror  Theorems \ref{main1} and \ref{main2}, respectively. The techniques in the $p=1$ case rely heavily on weighted energy estimates, the construction of suitable test functions adapted to the conformal geometry, and a careful application of the Laplacian comparison theorem.
All of our uniqueness results are ultimately derived from a general criterion formulated in terms of conformal geometry, see Theorems \ref{teo1} and \ref{teo2}. We reinterpret the original weighted heat equation as a classical heat equation on the conformally transformed manifold $(M,\tilde{g}, \tilde{\mu})$, where $\tilde{g} = \rho g$ and $d\tilde{\mu} = \rho d\mu$. In this setting, the weighted Laplacian becomes $\tilde{\Delta} = \frac{1}{\rho} \Delta$, and the heat equation takes the simplified form \eqref{introproblem}. 
\medskip

In the special case of model manifolds, see Section \ref{sec3} for a detailed definition, we demonstrate the sharpness of our results. In particular, we first provide in Proposition \ref{prop-nonuniq} sufficient conditions for the existence of infinitely many solutions, obtained via the construction of supersolutions to a related elliptic problem. Then, we construct explicit examples on model manifolds where our uniqueness theorems hold. Afterwards, we prove for such examples that, violating the assumptions on $\rho$ made in the uniqueness results, one can apply Proposition \ref{prop-nonuniq} showing that the Cauchy problem admits infinitely many nontrivial bounded solutions
These examples include standard model geometries such as the Euclidean space $\R^N$ and hyperbolic space $\mathbb H^N$, as well as more general Riemannian manifolds.

\

\noindent \textbf{Structure of the paper.} The main results of the paper are presented in Section \ref{sec1}, where the theorems are formulated distinguishing between the cases $p > 1$ and $p = 1$. Section \ref{sec3} provides auxiliary results and notational conventions, including the setting of weighted Riemannian manifolds, properties of model manifolds with a pole. In Section \ref{sec4}, the case $p > 1$ is treated in detail. We first state a general uniqueness theorem for the heat equation posed on a weighted Riemannian manifold. Then, by reducing problem \eqref{problema} to problem \eqref{heat}, we prove Theorems \ref{teogen}, \ref{main1} and \ref{main2}. A similar approach is used in Section \ref{sec5} to address the case $p = 1$, which requires additional geometric assumptions. Section \ref{sec-nonuniq} discusses a general non-uniqueness criteria in the special case of $M=\mathbb{M}^N_f$, i.e. a model manifold.
Finally, in Section \ref{sec-examples} we collect several illustrative examples that demonstrate the applicability of the main theorems and show that our assumptions on $\rho$ are sharp in the sense that, if they fail, than infinitely many solutions exist.

\bigskip

\
\

\section{Main results} \label{sec1}\setcounter{equation}{0}

Let $o$ be a reference point in $(M,g)$ and let 
$$
r(x)=\operatorname{dist}(x,o)\, \quad\text{for all}\,\,x\in M\, , 
$$
the Riemannian distance with respect to $g$. Throughout the paper we deal with classical solutions to equation \eqref{problema0}, we recall the definition.

\begin{definition}\label{defsole}
Let $u_0\in L^p_{\operatorname{loc}}(M)$, $f\in C(M\times(0,T])$, $p\ge1$. We say that a function $u$ is a classical solution to equation \eqref{problema0} if
\begin{itemize}
\item[(i)]  $u\in  C^{2,1}_{x,t}(M\times(0,T])$;
\item[(ii)] $\rho(x)\,\partial_tu(x,t)=\Delta u(x,t) + f(x,t) \,$\;\;  for all \;\; $(x,t)\in M\times(0,T]$;
\item[(iii)] $\lim_{t\to 0^+} \int_K|u(x,t)-u_0(x)|^p\,d\mu(x)=0$\;\; for any compact subset \;\; $K\subset M$.
\end{itemize}
Furthermore, we say that $u$ is a supersolution (subsolution) to equation \eqref{problema0}, if in $(ii)$ instead of $``="$ we have $``\geq"\; (``\leq")$\,.
\end{definition}

We now state the main uniqueness results distinguishing between the case of solutions belonging to a suitable weighted $L^p(M\times(0,T))$ space with $p>1$ and those in the suitable weighted $L^1(M\times(0,T))$. These two choices require different assumptions on the density $\rho$ and on the manifold $M$.
 
\subsection{Uniqueness results for $p>1$}

Let us first consider the cases where $\rho$ is bounded between two powerlike functions.

\begin{theorem}\label{main1}
Let $(M,g,\mu)$ be a complete, non-compact, weighted Riemannian manifold and $p>1$. Let $u$ be a solution to problem \eqref{problema} such that $\rho$ satisfies \eqref{compl} and, for some $0\le\theta< 2$
\begin{equation}\label{eq:rhobound-theta}
c_1(r+1)^{-\theta}\le\rho(x)\le c_2(r+1)^{-\theta}\quad\text{for all}\,\,\,x\in M\,,
\end{equation}
for some $c_1,c_2>0$.
Moreover, let $\psi$ be a positive, increasing and continuous function defined in $(0,+\infty)$ such that, for some $R_0>0$ big enough 
\begin{equation}\label{integral-phi-theta}
\int_{R_0}^\infty \frac { r^{1-\theta}}{\psi(r)}\,dr=+\infty\,.
\end{equation}
If $u\in L^p_{e^{-\psi}}(M\times(0,T))$,
then $$u\equiv 0\quad\text{in} \,\,\,M\times(0,T].$$
\end{theorem}

\begin{theorem}\label{main2}
Let $(M,g,\mu)$ be a complete, non-compact, weighted Riemannian manifold and $p>1$. Let $u$ be a solution to problem \eqref{problema} such that $\rho$ satisfies \eqref{compl} and 
\begin{equation}\label{eq:rhobound}
c_1(r+1)^{-2}\le\rho(x)\le c_2(r+1)^{-2}\quad\text{for all}\,\,x\in M\, , 
\end{equation}
for some $c_1,c_2>0$. Moreover, let $\psi$ be a positive, increasing and continuous function defined in $(0,+\infty)$ such that, for some $R_0>0$ big enough 
\begin{equation}\label{integral-phi}
\int_{R_0}^\infty \frac {\log r}{r\,\psi(r)}\,dr=+\infty\,.
\end{equation}
If $u\in L^p_{e^{-\psi}}(M\times(0,T))$,
then 
$$u\equiv 0\quad\text{in} \,\,\,M\times(0,T].$$
\end{theorem}

For more general choices of $\rho\in C^\infty(M)$, $\rho>0$ such that 
\begin{equation}\label{eq:rhobound-general}
\rho_1(r)\le\rho(x)\le\rho_2(r)\quad\text{for all}\,\,x\in M\,,
\end{equation} 
for some  $\rho_1,\rho_2:[0,+\infty)\to(0,+\infty)$, $\rho_1,\rho_2\in C^\infty([0,+\infty))$, we have the following

\begin{theorem}\label{teogen}
Let $(M,g,\mu)$ be a complete, non-compact, weighted Riemannian manifold and $p>1$. Let $u$ be a solution to problem \eqref{problema} such that $\rho$ satisfies \eqref{compl} and \eqref{eq:rhobound-general}. Moreover, let $\psi$ be a positive, increasing and continuous function defined in $(0,+\infty)$ such that, for some $R_0>0$ big enough 
\begin{equation}\label{integral-phi-gen}
\int_{R_0}^\infty \frac { \int_0^r\sqrt{\rho_1(s)}\,ds}{\psi(r)+\log(\rho_2(r))}\sqrt{\rho_1(r)}\,dr=+\infty\,;
\end{equation}
and \begin{equation}\label{monotonia}
\psi'(r)+\frac{\rho_2'(r)}{\rho_2(r)}\ge0 \quad\text{for every}\,\,\,r>R_0.\end{equation}
If $u\in L^p_{e^{-\psi}}(M\times(0,T))$;
then $$u\equiv 0\quad\text{in} \,\,\,M\times(0,T].$$
\end{theorem}

\begin{remark}
Observe that assumptions \eqref{eq:rhobound-theta} and \eqref{eq:rhobound} correspond to \eqref{eq:rhobound-general} with $\rho_1(r)=c_1(r+1)^{-\theta}$, $\rho_2(r)=c_2(r+1)^{-\theta}$ and to \eqref{eq:rhobound-general} with $\rho_1(r)=c_1(r+1)^{-2}$, $\rho_2(r)=c_2(r+1)^{-2}$ respectively. 
\end{remark}

\subsection{Uniqueness results for $p=1$}

To extend the class of uniqueness to the case $p=1$, one needs to require extra assumptions on the manifold $(M,g,\mu)$ and on the density function $\rho$. Firstly, we assume that $\rho\in C^\infty(M)$, $\rho>0$ and $\rho$ is radial, i.e.
\begin{equation}\label{eq:rhobound-general-p1}
\rho(x)\equiv \rho(r)\quad\text{for all}\,\,x\in M\,.
\end{equation}
This hypothesis is technical and strictly related to the method of proof. 
Similarly as before,  we first state the result in the special cases of $\rho$ being a power function which decays to zero at infinity.

\begin{theorem}\label{main1-p1}
Let $(M,g,\mu)$ be a complete, non-compact, weighted Riemannian manifold with a pole. Let $u$ be a solution to problem \eqref{problema} such that $\rho$ satisfies \eqref{compl} and, for some $0\le\theta< 2$,
\begin{equation}\label{eq:rhobound-theta-p1}
\rho(x)\equiv\rho(r)=c_1(r^2+1)^{-\frac{\theta}{2}}\quad\text{for all}\,\,x\in M\,,
\end{equation}
for some $c_1>0$.
Moreover, assume that
\begin{equation}\label{hp-laplaciano}
\Delta r\ge \frac{\theta(N-1)r}{2(1+r^2)}\quad \text{in } M\, .
\end{equation}
Let $\psi$ be a positive, increasing and continuous function defined in $(0,+\infty)$ such that, for some $R_0>0$ big enough 
\begin{equation}\label{integral-phi-theta-p1}
\int_{R_0}^\infty \frac { r^{1-\theta}}{\psi(r)}\,dr=+\infty\,.
\end{equation}
If $u\in L^1_{e^{-\psi}}(M\times(0,T))$, then $$u\equiv 0\quad\text{in} \,\,\,M\times(0,T].$$
\end{theorem}

\begin{theorem}\label{main2-p1}
Let $(M,g,\mu)$ be a complete, non-compact, weighted Riemannian manifold with a pole.  Let $u$ be a solution to problem \eqref{problema} such that $\rho$ satisfies \eqref{compl} and 
\begin{equation}\label{eq:rhobound-p1}
\rho(x)\equiv\rho(r)=\frac{c_2}{r^2+1}\quad\text{for all}\,\,x\in M\,,
\end{equation}
for some $c_2>0$. Moreover, assume that
\begin{equation}\label{hp-laplaciano-theta}
\Delta r\ge \frac{(N-1)r}{1+r^2}\quad \text{in } M\, . 
\end{equation}
Let $\psi$ be a positive, increasing and continuous function defined in $(0,+\infty)$ such that, for some $R_0>0$ big enough 
\begin{equation}\label{integral-phi-p1}
\int_{R_0}^\infty \frac {\log r}{r\,\psi(r)}\,dr=+\infty\,.
\end{equation}
if $u\in L^1_{e^{-\psi}}(M\times(0,T))$, then $$u\equiv 0\quad\text{in} \,\,\,M\times(0,T].$$
\end{theorem}
\medskip

We conclude the statements with the following Theorem for a more general density $\rho$.
\begin{theorem}\label{teogen-p1}
Let $(M,g,\mu)$ be a complete, non-compact, weighted Riemannian manifold with a pole. Let $u$ be a solution to problem \eqref{problema} such that $\rho$ satisfies \eqref{compl}, \eqref{eq:rhobound-general-p1} and
\begin{equation}\label{decayrho}
-\frac{\rho'(r)}{\rho(r)}=O(r)\quad \text{for all}\,\,\,r>1\, . 
\end{equation}
 Moreover, assume that  
\begin{equation}\label{hp-laplaciano-p1}
\Delta r\ge-\frac{N-1}{2}\frac{\rho'(r)}{\rho(r)}\quad\text{in }\,M\,  . 
\end{equation}
Let $\psi$ be a positive, increasing and continuous function defined in $(0,+\infty)$ such that, for some $R_0>0$ big enough 
\begin{equation}\label{integral-phi-gen-p1}
\int_{R_0}^\infty \frac { \int_0^r\sqrt{\rho(s)}\,ds}{\psi(r)+\log(\rho(r))}\sqrt{\rho(r)}\,dr=+\infty\,;
\end{equation}
and 
\begin{equation}\label{monotoniap1}
\psi'(r)+\frac{\rho'(r)}{\rho(r)}\geq 0 \quad \text{ for every } r>R_0\, . 
\end{equation}
If  $u\in L^1_{e^{-\psi}}(M\times(0,T))$,
then $$u\equiv 0\quad\text{in} \,\,\,M\times(0,T].$$
\end{theorem}

\noindent We now provide some examples of manifolds to which Theorems \ref{main1-p1}, \ref{main2-p1} and \ref{teogen-p1}, apply. The notation used in the following Remark relies on the definitions introduced in Section \ref{sec3}.
\begin{remark}
In Theorems \ref{main1-p1}, \ref{main2-p1} and \ref{teogen-p1} one can choose:
\begin{itemize}
\item Any model manifold $\mathbb{M}^N_f$ with $f(r)\geq \frac{1}{\sqrt{\rho(r)}}$ for all $r\geq 0$. In this case, 
$$
\Delta r=(N-1)\frac{f'(r)}{f(r)}
$$
therefore  \eqref{hp-laplaciano}, \eqref{hp-laplaciano-theta} and \eqref{hp-laplaciano-p1} are satisfied. 
\item A generic Riemannian manifold with a pole such that 
$$
\mathrm{sec}_{\mathrm{rad}}(x)\leq -\dfrac{h''(r(x))}{h(r(x))}\, , \quad \text{ for all } x\in M\setminus\{o\}\,  ,
$$
holds, where $h(r)\geq \frac{1}{\sqrt{\rho(r)}}$ for all $r\geq 0$. In this case  \eqref{hp-laplaciano}, \eqref{hp-laplaciano-theta} and \eqref{hp-laplaciano-p1} follow from the Laplacian comparison, see Theorem \ref{laplacian_comp}. 
\end{itemize}
As a remarkable example we point out that in the Theorems \ref{main1-p1} and \ref{main2-p1}, the manifold can be chosen to be  the hyperbolic space $\mathbb{H}^N$. Indeed, it is a complete, non-compact Riemannian manifold with a pole and 
$$
\Delta r=(N-1)\coth(r) \geq (N-1)\frac{r}{1+r^2}\geq \frac{\theta(N-1)}{2}\frac{r}{1+r^2}\, ,
$$
therefore \eqref{hp-laplaciano} and \eqref{hp-laplaciano-theta} are satisfied. For a proper choice of $\rho$, we can chose $M=\mathbb H^N$ also in Theorem \ref{teogen-p1}.
\end{remark}

\section{Notations from Riemannian Geometry}\label{sec3}\setcounter{equation}{0}

In this section we recall and fix some basic notations and useful definitions from Riemannian Geometry that we will use in the paper (we refer e.g. to \cite{AMR,Chavel,Grigbook,Petersen,PRS}).

\

\noindent \textbf{Basic Notions.} Throughout the paper, $M$ will denote a smooth $N-$dimensional, complete, noncompact, Riemannian manifold without boundary and of infinite volume, endowed with a smooth Riemannian metric $g=\lbrace g_{ij}\rbrace$. As usual, we denote by $\mathrm{div}(X)$ the \emph{divergence} of a smooth vector field $X$ on $M$, that is,  in local coordinates
\begin{equation*}
\mathrm{div} (X)=\dfrac{1}{\sqrt{\vert g\vert}}\sum_{i=1}^{N} \partial_i \left( \sqrt{\vert g\vert}X ^i\right)\, ,
\end{equation*}
where $\vert g\vert=\mathrm{det}(g_{ij})(\geq 0)$. We also denote by $\nabla$ and $\Delta$ the \emph{gradient} and the \emph{Laplace-Beltrami operator} on $M$ with respect to $g$, that is, in local coordinates,
$$
(\nabla u)^i=\sum_{j=1}^{N}g^{ij}\partial_j u \, ,
$$
and
$$
\Delta u=\mathrm{div}(\nabla u)=\dfrac{1}{\sqrt{\vert g\vert}}\sum_{i,j=1}^{N} \partial_i \left( \sqrt{\vert g\vert} g^{ij}\partial_j u\right)\, ,
$$
for any $C^2-$function $u:M\rightarrow\mathbb{R}$, where $\lbrace g^{ij}\rbrace$ denotes the inverse of the metric tensor $g$. Moreover, let $\mu$ denotes the Riemannian volume form on $M$ with respect to $g$, that is, in local coordinates
$$
d\mu=\sqrt{|g|}dx^1\dots dx^N\, . 
$$
Due to the divergence structure of the Laplace-Beltrami operator we have the \emph{integration by parts formula} 
$$
\int_M v\Delta u\, d\mu=- \int_M g(\nabla u,\nabla v)\, d\mu\, ,
$$
for any $C^2-$functions $u,v:M\rightarrow\mathbb{R}$, with either $u$ or $v$ compactly supported.

\

\noindent \textbf{Weighted Manifolds.} We recall that a \emph{weighted manifold} (or manifold with density, or smooth metric measure space) is a triad $(M,g,\mu_\omega)$ where $(M,g)$ is a Riemannian manifold with volume form $\mu$ and 
\begin{equation}\label{changemeasure}
d\mu_\omega =e^{-\omega} d\mu\, ,
\end{equation}
where $\omega:M\rightarrow\mathbb{R}$ is a smooth function. The definition of gradient on a weighted Riemannian manifold $(M,g,\mu_ \omega)$ is the same as on $(M,g)$, but the definition of the divergence changes. Given a smooth vector field $X$ on $M$ we define the \emph{weighted divergence} by
$$
\mathrm{div}_\omega(X)=e^\omega\mathrm{div}\left(e^{-\omega}X\right)\, ,
$$
and the \emph{weighted Laplacian} of a $C^2-$function $u:M\rightarrow\mathbb{R}$ by 
$$
\Delta_\omega u=\mathrm{div}_\omega(\nabla u)=\Delta u-g(\nabla \omega,\nabla u)\, .
$$ 
It is easy to see that the following \emph{weighted integration by parts formula} 
\begin{equation}\label{int_parts_f}
\int_Mu\,\Delta_{\omega} v\,d\mu_\omega=-\int_M g(\nabla u,\nabla v)\,d\mu_\omega\,,
\end{equation}
holds, provided $u,v:M\rightarrow\mathbb{R}$ are $C^2-$functions, with either $u$ or $v$ compactly supported.

Observe that, for any $\psi\in C^2(\R)$, $v\in C^2(M)$, we have
\begin{equation}\label{deltacomposition}
\Delta_{\omega}[\psi(v)]=\psi'(v)\Delta_{\omega} v+\psi''(v)|\nabla v|^2.
\end{equation}

\

\noindent \textbf{Conformal metric.} Given a Riemannian manifold $(M,g)$, we consider the function $\rho\in C^\infty(M)$, $\rho>0$; then define the \emph{conformal metric}
\begin{equation}\label{conformalchange}
\tilde g=\rho\,g\, .
\end{equation}
The natural assumption that guarantees that the Riemannian manifold $(M,\tilde{g})$ is complete is \eqref{compl}, i.e. 
\begin{equation*}
\int_1^\infty\sqrt{\rho(\gamma(s))}\, ds=+\infty\, ,\quad \text{ for every divergent ray $\gamma$}\,.
\end{equation*}
We recall that a minimizing geodesic $\gamma:[0,+\infty)\rightarrow M$ is called a \emph{ray} and it is called \emph{divergent} if eventually it leaves every compact subset of $M$. Indeed, from Hopf-Rinow Theorem we know that $(M,\tilde{g})$ is complete if and only if every ray has infinite $\tilde{g}-$length; consider a ray $\gamma:[0,+\infty)\rightarrow M$ and we parametrize it by $g-$arclength, i.e. 
$$
|\dot{\gamma}(s)|_g=1\, ,
$$ 
then 
$$
|\dot{\gamma}(s)|_{\tilde{g}}=\sqrt{\tilde{g}(\dot{\gamma}(s),\dot{\gamma}(s))}=\sqrt{\rho(\gamma(s))}\, .
$$ 
In particular,
$$
L_{\tilde{g}}(\gamma)=\int_0^\infty|\dot{\gamma}(s)|_{\tilde{g}}\, ds=\int_0^\infty \sqrt{\rho(\gamma(s))}\, ds\geq \int_1^\infty \sqrt{\rho(\gamma(s))}\, ds=+\infty\, ,
$$
from \eqref{compl}.
Moreover, we consider the weighted Riemannian manifold obtained by setting in \eqref{changemeasure} 
$\omega =\log\frac{1}{\rho}$, i.e. 
\begin{equation}\label{weighted_mess}
d\tilde\mu_\omega=\rho\,d\mu\,.
\end{equation}
We denote the triplet $(M,\tilde{g},\tilde{\mu}_\omega)$, with a slight abuse of notation, as $$(M,\tilde{g},\tilde{\mu}_\rho)$$ so that the dependence on $\rho$ is evident. 
It is not difficult to show that (see e.g. \cite[Chapter 3]{Grigbook}) the weighted Laplacian with respect to $\tilde{g}$ is
\begin{equation}\label{deltaweighted}
\tilde\Delta:=\tilde\Delta_{\omega}=\frac1{\rho}\,\Delta\, .
\end{equation}
By \eqref{deltaweighted}, problem \eqref{problema} can be written as
\begin{equation}\label{problemachanged}
\begin{cases}
\partial_t u=\tilde \Delta u \quad  &\text{in}\,\,\,M\times(0,T],\\
u=0\quad &\text{on}\,\,\,M\times\{0\}\,,
\end{cases}
\end{equation}
where $\tilde g$ and $\tilde \mu$ are given by \eqref{conformalchange} and \eqref{weighted_mess}.

\begin{remark}\label{regularityrho}
Observe that the assumption on the smoothness of $\rho$ is necessary because $\rho$ is involved into the conformal change of the metric $g$ and the measure $\mu$ into $\tilde g$ and $\tilde \mu$ (see \eqref{conformalchange}). Although, when dealing only with uniqueness of solutions to problem \eqref{problema}, one does not need to require more than continuity on the function $\rho$.
\end{remark}

\

\noindent \textbf{Manifolds with pole and Laplacian comparison.} Let $(M,g)$ be a connected, complete Riemannian manifold of dimension $N\geq 2$. For every $x\in M$, let $r(x)$ be the distance function from a reference point $o\in M$, it is well-known that $r(x)$ is a Lipschitz function on $M$ and smooth on $M\setminus \left(\{o\}\cup \mathrm{cut}\{o\}\right)$ with $|\nabla r|\leq 1$, where $\mathrm{cut}\{o\}$ denotes the cut-locus of $o$. We recall that $o$ is called a \emph{pole} if $\mathrm{cut}\{o\}=\emptyset$.

\begin{definition}\label{def_modello}
An $N-$dimensional model manifold with warping function $\varphi$ is a Riemannian manifold $\mathbb{M}^N_{f}$ diffeomorphic to $\mathbb{R}^N$ and endowed with a rotationally symmetric Riemannian metric. The model $\mathbb{M}^N_{f}$ is realized as the quotient space $[0,+\infty)\times\mathbb{S}^{N-1}/\backsim$, where $\backsim$ identifies $\lbrace 0\rbrace\times\mathbb{S}^{N-1}$ with the pole $o$ of the space, and the Riemannian metric has the expression
$$
g=dr^2+f^2(r)d\theta^2\, 
$$
where $r(x)=\operatorname{dist}(o,x)$,  $d\theta^2$ denotes the standard metric of $\mathbb{S}^{N-1}$ and $f:[0,+\infty)\rightarrow[0,+\infty)$ is a smooth function satisfying the following conditions: 
\begin{itemize}
\item $f(r)>0$ for all $r>0$,
\item $f'(0)=1$, 
\item $f^{(2k)}(0)=0$ for every $k\geq 0$.
\end{itemize}
\end{definition}
Observe that the Euclidean space $\mathbb{R}^N$ and the (generalized) hyperbolic space $\mathbb{H}^N$ of constant curvature $-c<0$ can be seen as model manifolds (take $f(r)=r$ and $f(r)=\frac{\sinh(\sqrt{c}r)}{\sqrt{c}}$, respectively). Moreover, due to their structure the Laplace-Beltrami operator on $\mathbb{M}^N_{\varphi}$ can be written in the following way
\begin{equation}\label{Laplaciano_modelli}
\Delta = \dfrac{\partial^2}{\partial^2 r} + (N-1)\dfrac{f'}{f}\dfrac{\partial}{\partial r} + \dfrac{1}{f^2}\Delta_\theta\, , 
\end{equation}
where $\Delta_\theta$ denotes the Laplace-Beltrami operator on $\mathbb{S}^{N-1}$. In particular, 
\begin{equation}\label{Laplaciano_dist}
\Delta r= (N-1)\dfrac{f'}{f}\, .
\end{equation}
Moreover, being the volume form on $\mathbb{M}^N_{f}$ given by
\begin{equation}\label{eq84}
d\mu=f^{N-1}(r)\, dr\, d\theta\, ,
\end{equation}
one has the following: 
\begin{equation}\label{volume_bolle}
V(o,r)=c_n\int_{0}^{r}f^{N-1}(t)\, dt \, ,
\end{equation}
where, $c_N$ is the area of the $(N-1)-$dimensional unit sphere and $o$ is the pole of $\mathbb{M}^N_{f}$.

In addition, the \emph{radial sectional curvature} of $\mathbb{M}^N_{f}$ is the following:
$$
\mathrm{sec}_{\mathrm{rad}}(x):=\mathrm{sec}\left(\nabla r,X\right)=-\frac{f''(r(x))}{f(r(x))}\, , 
$$
where $\mathrm{sec}_{\mathrm{rad}}(x)$ denotes the sectional curvature at a point $x\in M$ of a $2-$plane determined by $\nabla r$ and $X\in\nabla r(x)^{\perp}$. Furthermore, the \emph{radial Ricci curvature} of $\mathbb{M}^N_{f}$ is
$$
\mathrm{Ric}(\nabla r(x),\nabla r(x))=-(N-1)\frac{f''(r(x))}{f(r(x))}\, .
$$
Having the previous expressions in mind we can recall the Laplacian comparison from below and above. We define the following class of functions 
\begin{equation}\label{A}
\mathcal{A}=\left\{ h\in C^{\infty}(0,\infty)\cap C^1([0,\infty)) \, : \, h(0)=0 \, , \, \, h'(0)=1 \, ,\, \,  h>0 \,  \text{ in }\,  (0,\infty)\right\}\, .
\end{equation}

\begin{theorem}[Laplacian comparisons]\label{laplacian_comp}
Let $(M,g)$ be a complete Riemannian manifold of dimension $N$ with a pole $o$. Consider a function $h\in\mathcal{A}$. 
\begin{itemize}
\item[$i)$] Assume that the radial sectional curvature satisfy
$$
\mathrm{sec}_{\mathrm{rad}}(x)\leq -\dfrac{h''(r(x))}{h(r(x))}\, , \quad \text{ for all } x\in M\setminus\{o\}\,  ,
$$
then
$$
\Delta r(x)\geq (N-1)\dfrac{h'(r(x))}{h(r(x))} \,  , \quad \text{ for every } x\in M\setminus\{o\}\, .
$$
\item[$ii)$] Assume that the radial Ricci curvature satisfy
$$
\mathrm{Ric}(\nabla r(x),\nabla r(x))\geq -(N-1)\dfrac{h''(r(x))}{h(r(x))}\, , \quad \text{ for all } x\in M\setminus\{o\}\,  ,
$$
then
$$
\Delta r(x)\leq (N-1)\dfrac{h'(r(x))}{h(r(x))} \,  , \quad \text{ for every } x\in M\setminus\{o\}\, . 
$$
\end{itemize}
\end{theorem}

We recall that for a manifold with a pole we can consider \emph{geodesic polar coordinates} with respect to the pole $o$, i.e. for any $x\in M\setminus\{o\}$ we have a polar radius $r=d_g(o,x)$ and a polar angle $\theta\in\mathbb{S}^{N-1}$ such that the shortest geodesics from $o$ to $x$ start at $o$ with the direction $\theta$ in the tangent space $T_oM$. Since we can identify $T_oM$ with $\mathbb{R}^n$, then $\theta$ can be regarded as a point on $\mathbb{S}^{n-1}$. In particular, $M\setminus\{o\}$ is diffeomorphic to a star-like region of $\mathbb{R}_{+}\times\mathbb{S}^{n-1}$. Moreover, the Riemannian metric $g$ in  $M\setminus\{o\}$ has, in the polar coordinates, the form
$$
g=dr^2+A(r,\theta)d\theta^2\, 
$$
where $d\theta^2$ denotes the standard metric in $\mathbb{S}^{n-1}$ and $A=\{A_{ij}\}$ denotes a positive definite matrix. It is easy to see that the Laplace-Beltrami operator in the polar coordinates has the form:
\begin{equation}\label{lapl_polar}
\Delta=\frac{\partial^2}{\partial r^2} + m(r,\theta)\frac{\partial}{\partial r} +  \Delta_{\theta}\, ,
\end{equation}
where $m(r,\theta)=\frac{\partial}{\partial r}\left(\log \sqrt{|A|}\right)$, $|A|=\mathrm{det}(A)$ and $\Delta_\theta$ is the Laplace-Beltrami operator on the submanifold $\partial B(o,r)$. In particular we have that 
\begin{equation}\label{m}
\Delta r=m(r,\theta)\, . 
\end{equation}

\section{Proof of Theorems  \ref{main1}, \ref{main2} and \ref{teogen}}\label{sec4}\setcounter{equation}{0}

In order to prove Theorems \ref{main1}, \ref{main2} and \ref{teogen}, we state a general uniqueness theorem for the Cauchy problem of the heat equation introduced in \eqref{problemachanged} where $M$ is a weighted, complete Riemannian manifold. The goal will be, by means of a suitable conformal change of the metric defined in \eqref{conformalchange}, to reformulate problem \eqref{problema} into problem \eqref{problemachanged} where the operator is given in \eqref{deltaweighted} and the measure is \eqref{weighted_mess}.

\subsection{Uniqueness for the heat equation on a weighted manifold with $p>1$}

The following theorem is formulated for any weighted, geodetically complete Riemannian manifold. In this subsection, the weight $e^{-\omega}$ does not depend on $\rho$ (where we have used the notation introduced in Section \ref{sec3}). For this reason, the following theorem is of independent interest.

\begin{theorem}\label{teo1}
Let $(M, g,\mu_\omega)$ be a weighted Riemannian manifold and $p>1$. Moreover, let $\varphi:M\to(0,+\infty)$ be an increasing and continuous function such that, for some $x_0\in M$ and $R_0>0$,
\begin{equation}\label{integral-varphi}
\int_{R_0}^\infty \frac{ r}{\varphi(r)}\,d r=+\infty\,.
\end{equation}
Let $u$ be a solution to problem \eqref{problemachanged} such that
\begin{equation}\label{mainassumptiontilde}
\int_0^T\int_{B_R(x_0)}|u(x,t)|^p\,d\mu_\omega(x)dt\le e^{\varphi(R)}\quad\text{for all}\,\,R>R_0\,.
\end{equation}
Then 
$$
u\equiv 0 \quad \text{ in } \quad M\times(0,T]\, .
$$
\end{theorem}

The proof of Theorem \ref{teo1} can be found in \cite[Theorem 2.2]{Punzo2015} where uniqueness is established for the heat equation of the Laplace-Beltrami operator posed on a geodetically complete Riemannian manifold. In that context, the manifold is not weighted but the proof can be adapted also to the weighted setting and the results reads the same. For this reasons we omit the proof.
\bigskip

\noindent The proofs of Theorems \ref{main1}, \ref{main2} and \ref{teogen} are based on a conformal change of the metric: since $\rho\in C^\infty(\R)$, $\rho>0$, it is possible to define the conformal metric $\tilde g$ and the measure $\tilde \mu_\omega$, with $\omega=\log\frac{1}{\rho}$ such that
\begin{equation}\label{eq75}
\tilde g=\rho\,g,\quad\quad d\tilde\mu_\omega:=d\tilde\mu=\rho\,d\mu\,.
\end{equation}
We therefore consider the conformal and weighted Riemannian manifold $(M,\tilde{g},\tilde{\mu}).$ Observe that, since \eqref{compl} is in force and $(M,g)$ is complete, we also have that $(M,\tilde{g})$ is complete. 
Furthermore, we recall that, due to \eqref{deltaweighted}, problem \eqref{problema} is in turn
\begin{equation}\label{problema:support}
\begin{cases}
\partial_tu=\tilde\Delta u\quad&\text{in}\,\,\,M\times(0,T]\\
u=0\quad &\text{on}\,\,\,M\times\{0\}\,;
\end{cases}
\end{equation}
therefore, if $u$ is a solution to \eqref{problema}, then it is also a solution to \eqref{problema:support}. For this reason, the goal of the following proofs is to find a proper $\varphi$ which allows us to apply Theorem \ref{teo1}.

\

\noindent \textbf{Proof of Theorem \ref{main1}.} 
 Let 
$$R_2:=\min\{r\ge R_0,\,\,:\,\,\psi(r)+\log(\rho(r))\ge0\}.$$ 
We define the function $\varphi:\R^+\to\R^+$
\begin{equation}\label{eq:phi}
\varphi(\tilde r)=\begin{cases}
\psi(R_2)+\log(c_2(1+R_2)^{-\theta}) \quad &\text{for}\,\,\,\tilde r\leq \tilde R_2 \\
\psi(r(\tilde r))+\log(c_2(1+r(\tilde r))^{-\theta}) \quad &\text{for}\,\,\,\tilde r>\tilde R_2
\end{cases}
\end{equation}
where
\begin{equation}\label{eq:tilde-r}
\sqrt{c_1}\int_{0}^{r(x)}(1+s)^{-\theta/2}\, ds\le\tilde r(r)\le\sqrt{c_2}\int_{0}^{r(x)}(1+s)^{-\theta/2}\, ds ;
\end{equation} 
and
\begin{equation}\label{eq:r}
\sqrt{c_2}\int_{0}^{\tilde r(x)}(1+s)^{\theta/2}\, ds\le  r(\tilde r)\le\sqrt{c_1}\int_{0}^{\tilde r(x)}(1+s)^{\theta/2}\, ds\, ;
\end{equation} 
with
$$
\tilde R_2=\sqrt{c_2}\int_{0}^{R_2}(1+s)^{-\theta/2}\, ds\,.
$$
Observe that $\varphi$ is continuous and increasing in the variable $\tilde r$ since $\psi$ and the function $\log$ are increasing as well as $r(\tilde r)$ since it is an integral function. In addition, $\varphi$ is non-negative. Furthermore, the solution $u$ for problem \eqref{problema:support} satisfies \eqref{mainassumptiontilde} with this choice of $\varphi$. In fact, since $u\in L^p_{e^{-\psi}}(M\times(0,T),d\mu)$ and due to  \eqref{eq:rhobound}, \eqref{eq75} and \eqref{eq:phi}, we have that, for every $R>\tilde R_2$ big enough,
$$
\begin{aligned}
\infty>C&\ge\int_0^T\int_M |u(x,t)|^p\,e^{-\psi}\,d\mu(x)dt\\
&=
\int_0^T\int_{B_{R_2}} |u(x,t)|^p\,e^{-\psi}\,d\mu(x)dt+\int_0^T\int_{M\setminus B_{\tilde R_2}} |u(x,t)|^p\,e^{-\varphi+\log(\rho_2)}\,\frac{1}{\rho}d\tilde\mu(x)dt\\
&=
\int_0^T\int_{B_{R_2}} |u(x,t)|^p\,e^{-\psi}\,d\mu(x)dt+\int_0^T\int_{M\setminus B_{\tilde R_2}} |u(x,t)|^p\,e^{-\varphi}\,\frac{\rho_2}{\rho}d\tilde\mu(x)dt\\
&\ge e^{-\psi(R_2)}\int_0^T\int_{B_{R_2}} |u(x,t)|^p\,d\mu(x)dt+\int_0^T\int_{M\setminus B_{\tilde R_2}} |u(x,t)|^p\,e^{-\varphi}\,d\tilde\mu(x)dt\\
&\ge e^{-\varphi(\tilde R_2)}c_2(1+R_2)^{-\theta}\int_0^T\int_{B_{\tilde R_2}} |u(x,t)|^p\,\frac{1}{\rho}d\tilde\mu(x)dt+\int_0^T\int_{M\setminus B_{\tilde R_2}} |u(x,t)|^p\,e^{-\varphi}\,d\tilde\mu(x)dt\\
&\ge e^{-\varphi(\tilde R_2)}c_2(1+R_2)^{-\theta}\int_0^T\int_{B_{\tilde R_2}} |u(x,t)|^p\frac{(1+r)^{\theta}}{c_2}\,d\tilde\mu(x)dt+\int_0^T\int_{M\setminus B_{\tilde R_2}} |u(x,t)|^p\,e^{-\varphi}\,d\tilde\mu(x)dt\\
&\ge e^{-\varphi(\tilde R_2)}(1+R_2)^{-\theta}\int_0^T\int_{B_{\tilde R_2}} |u(x,t)|^p\,d\tilde\mu(x)dt+\int_0^T\int_{M\setminus B_{\tilde R_2}} |u(x,t)|^p\,e^{-\varphi}\,d\tilde\mu(x)dt\\
&=(1+R_2)^{-\theta}\int_0^T\int_{B_{\tilde R_2}} |u(x,t)|^pe^{-\varphi}\,d\tilde\mu(x)dt+\int_0^T\int_{M\setminus B_{\tilde R_2}} |u(x,t)|^p\,e^{-\varphi}\,d\tilde\mu(x)dt\\
&\ge (1+R_2)^{-\theta}\int_0^T\int_{B_{\tilde R_2}} |u(x,t)|^pe^{-\varphi}\,d\tilde\mu(x)dt+ \int_0^T\int_{B_R\setminus B_{\tilde R_2}} |u(x,t)|^pe^{-\varphi}\,\,d\tilde\mu(x)dt\\
&\ge \tilde C\int_0^T\int_{B_R} |u(x,t)|^pe^{-\varphi}\,\,d\tilde\mu(x)dt \ge \tilde C e^{-\varphi(R)} \int_0^T\int_{B_R} |u(x,t)|^p\,\,d\tilde\mu(x)dt\,,
\end{aligned}
$$
for some $\tilde C=\tilde C(R_2,\theta)$. Thus \eqref{mainassumptiontilde} holds, in fact the latter yields, for some $C>0$
$$
\int_0^T\int_{B_R} |u(x,t)|^p\,\,d\tilde\mu(x)dt\le C\,e^{-\varphi(R)}\quad \text{for every $R>\tilde R_2$}\,.
$$
We are left to show that $\varphi$ defined in \eqref{eq:phi} satisfies \eqref{integral-varphi}. Due to \eqref{eq:rhobound-theta} and \eqref{eq:tilde-r} we have
\begin{equation}\label{r-tilda-theta}
\begin{aligned}
\frac{2c_1}{2-\theta}(1+r(x))^{(2-\theta)/2}
\le\tilde r(x)
\le\frac{2c_2}{2-\theta}(1+r(x))^{(2-\theta)/2}\,;
\end{aligned}
\end{equation}
 and consequently,
\begin{equation}\label{oppositetheta}
\left(\frac{2-\theta}{2c_2}\,\tilde r\right)^{2/(2-\theta)}-1\le r \le \left(\frac{2-\theta}{2c_1}\,\tilde r\right)^{2/(2-\theta)}-1\,.
\end{equation}
We consider \eqref{integral-phi-theta} and we perform the change of variable in \eqref{eq:tilde-r}, therefore we get
\begin{equation}\label{eq77}
+\infty=\int_{R_0}^\infty\frac{r^{1-\theta}}{\psi(r)}\,dr
\le\int_{\tilde R_0}^\infty\frac{(r(\tilde r))^{1-\theta}}{\psi(r(\tilde r))}\,\frac{d\tilde r}{\sqrt{\rho_1}(1+r(\tilde r))^{-\theta/2}}\,,
\end{equation}
with 
$$
\tilde R_0=\int_{0}^{R_0}\sqrt{c_1}(1+s)^{-\theta/2}\,ds\,.
$$
Let $\infty>\tilde R_1\ge\max\left\{\tilde R_0,\,\frac2{2-\theta}c_2^{2-\theta/2}\right\}$. Then observe that, since the integrand function is continuous, we have
$$
\int_{\tilde R_0}^{\tilde R_1}\frac{(r(\tilde r))^{1-\theta}}{\psi(r(\tilde r))}\,\frac{d\tilde r}{\sqrt{c_1}(1+r(\tilde r))^{-\theta/2}}<\infty\,,
$$
and therefore, due to \eqref{eq77},
\begin{equation}\label{eq01}
+\infty=\int_{\tilde R_1}^\infty\frac{(r(\tilde r))^{1-\theta}}{\psi(r(\tilde r))}\,\frac{d\tilde r}{\sqrt{c_1}(1+r(\tilde r))^{-\theta/2}}\,.
\end{equation}
By \eqref{eq:rhobound}, \eqref{eq:phi}, \eqref{r-tilda-theta} and \eqref{oppositetheta}, inequality \eqref{eq01}, for some constant $C=C(\tilde R_0,c_1,\theta)>0$, yields
$$
\begin{aligned}
+\infty&=\int_{\tilde R_1}^\infty\frac{(r(\tilde r))^{1-\theta}}{\psi(r(\tilde r))}\,\frac{d\tilde r}{\sqrt{c_1}(1+r(\tilde r))^{-\theta/2}}\\
&\le \frac{1}{\sqrt{c_1}} \int_{\tilde R_1}^\infty \frac{(r(\tilde r))^{1-\theta}}{[\varphi(\tilde r)-\log(c_2(1+r(\tilde r))^{-\theta})]}{(1+r(\tilde r))^{\theta/2}}\,d\tilde r\\
&\le\frac{C_1}{\sqrt{c_1}} \int_{\tilde R_1}^\infty \frac{(1+r(\tilde r))^{1-\theta/2}}{\varphi(\tilde r)+\theta\log\left((1+r(\tilde r))\right)-\log(c_2)}\,d\tilde r\\
&\le\frac{C_2(2-\theta)}{2c_1^{3/2}} \int_{\tilde R_1}^\infty \frac{\tilde r}{\varphi(\tilde r)+\frac{2\theta}{2-\theta}\log\left(\frac{2-\theta}{2}\,\frac{\tilde r}{c_2}\right)-\log c_2}\,d\tilde r\\
&\le C \int_{\tilde R_1}^\infty \frac{\tilde r}{\varphi(\tilde r)}\,d\tilde r\,;
\end{aligned}
$$
where we have used that 
$$
\frac{2\theta}{2-\theta}\log\left(\frac{2-\theta}{2}\,\frac{\tilde r}{c_2}\right)\ge\log c_2\quad \text{ for all }\,\,\tilde r\ge \tilde R_1.
$$
 Thus we have showed that $\varphi$ satisfies \eqref{integral-varphi} for all $\tilde r> \tilde R_1$. Finally, by picking $R_0$ in Theorem \ref{teo1} to be the maximum between $\tilde R_1$ and $\tilde R_2$, the thesis follows.

\

\noindent \textbf{Proof of Theorem \ref{main2}.} 
We define $\varphi$ as in \eqref{eq:phi} with $\theta=2$. The same argument used in the proof of Theorem \ref{main1}, shows that $u$ satisfies \eqref{mainassumptiontilde} with this choice of $\varphi$. We are left to show that $\varphi$ satisfies \eqref{integral-varphi}. Due to \eqref{eq:rhobound}, 
\begin{equation}\label{r-tilda-2}
\begin{aligned}
c_1\log(1+r(x))
\le\tilde r(x)\le c_2\log(1+r(x))\,;
\end{aligned}
\end{equation}
 and consequently, we have
\begin{equation}\label{opposite2}
e^{\frac{\tilde r}{c_2}}-1\le r \le e^{\frac{\tilde r}{c_1}}-1\,.
\end{equation}
We consider \eqref{integral-phi} and we perform the change of variable in \eqref{eq:tilde-r}, therefore we get
\begin{equation}\label{eq76}
+\infty=\int_{R_0}^\infty\frac{\log r}{r\,\psi(r)}\,dr\le\int_{\tilde R_0}^\infty\frac{\log (r(\tilde r))}{r(\tilde r)\,\psi(r(\tilde r))}\frac{1+r(\tilde r)}{\sqrt{c_1}}\,d\tilde r\,,
\end{equation}
with 
$$
\tilde R_0=\int_{0}^{R_0}\frac{\sqrt{c_1}}{1+s}\,ds\,.
$$
By \eqref{eq:rhobound}, \eqref{eq:phi}, \eqref{r-tilda-2} and \eqref{opposite2}, equality \eqref{eq01}, for some constant $C=C(R_0)>\frac1{R_0}+1$, yields
$$
\begin{aligned}
+\infty&=\int_{\tilde R_0}^\infty\frac{\log (r( \tilde r))}{r(\tilde r)\,\psi(r(\tilde r))}\,\frac{1+r(\tilde r)}{\sqrt{c_1}}\,d\tilde r\\
&= \frac{1}{\sqrt{c_1}} \int_{\tilde R_0}^\infty \frac{\log\left(e^{\frac{\tilde r}{c_1}}\right)}{[\varphi(\tilde r)-\log(c_1(1+r(\tilde r))^{-2})]}\frac{1+r(\tilde r)}{r(\tilde r)}\,d\tilde r\\
&\le \frac{1}{C c_1}  \int_{\tilde R_0}^\infty \frac{\tilde r}{\varphi(\tilde r)+2\log\left(1+r(\tilde r)\right)-\log(c_2)}\,d\tilde r\\
&\le \frac{1}{C c_1} \int_{\tilde R_0}^\infty \frac{\tilde r}{\varphi(\tilde r)+2\frac{\tilde r}{c_2}-\log c_2}\,d\tilde r\\
&\le \frac{1}{C c_1} \int_{\tilde R_0}^\infty \frac{\tilde r}{\varphi(\tilde r)}\,d\tilde r\,;
\end{aligned}
$$

where we have used that $(1+r)/r\le C$ for all $r\ge R_0$. This completes the proof.

\

\noindent \textbf{Proof of Theorem \ref{teogen}.}  Let 
$$R_2:=\min\{r\ge R_0,\,\,:\,\,\psi(r)+\log(\rho_2(r))\ge0\}.$$ 
We define the function $\varphi:\R^+\to\R^+$
\begin{equation}\label{eq:phi-gen}
\varphi(\tilde r)=\begin{cases}
\psi(R_2)+\log(\rho_2(R_2)) \quad &\text{for}\,\,\,\tilde r\leq \tilde R_2 \\
\psi(r)+\log(\rho_2(r)) \quad &\text{for}\,\,\,\tilde r>\tilde R_2
\end{cases}
\end{equation}
with
$$
\tilde R_2=\int_{0}^{R_2}\rho_2(s)\, ds\,.
$$
Moreover, we have the relations
\begin{equation}\label{eq:tilde-r-gen}
\int_{0}^{r(x)}\sqrt{\rho_1(s)}\, ds\le\tilde r(r)\le\int_{0}^{r(x)}\sqrt{\rho_2(s)}\, ds ;
\end{equation} 
and
\begin{equation}\label{eq:r-gen}
\int_{0}^{\tilde r(x)}\left(\sqrt{\rho_2(s)}\right)^{-1}\, ds\le  r(\tilde r)\le\int_{0}^{\tilde r(x)}\left(\sqrt{\rho_1(s)}\right)^{-1}\, ds\, ,
\end{equation} 
see e.g. \cite[Section 3.11]{Grigbook}. Observe that $\varphi$ is continuous and increasing in the variable $\tilde r$ due to the assumptions on $\psi$ and $\rho$, see \eqref{monotonia}.
In addition, $\varphi$ is non-negative. Furthermore, arguing as in the proof of Theorem \ref{main1}, we can infer that the solution $u$ to problem \eqref{problema:support} satisfies \eqref{mainassumptiontilde} with this choice of $\varphi$. We are left to show that $\varphi$ defined in \eqref{eq:phi-gen} satisfies \eqref{integral-varphi}.
We consider \eqref{integral-phi-gen}, by using \eqref{eq:phi-gen} and \eqref{eq:tilde-r-gen} we get
\begin{equation}\label{eqA1}
\begin{aligned}
+\infty&=\int_{R_0}^\infty \frac { \int_0^r\sqrt{\rho_1(s)}\,ds}{\psi(r)+\log(\rho_2(r))}\sqrt{\rho_1(r)}\,dr\\
&=\int_{R_0}^{R_2} \frac { \int_0^r\sqrt{\rho_1(s)}\,ds}{\psi(r)+\log(\rho_2(r))}\sqrt{\rho_1(r)}\,dr+\int_{R_2}^\infty \frac { \int_0^r\sqrt{\rho_1(s)}\,ds}{\psi(r)+\log(\rho_2(r))}\sqrt{\rho_1(r)}\,dr\\
&\le C+\int_{\tilde R_2}^\infty\frac{\tilde r}{\varphi(\tilde r)}\,d\tilde r\,.
\end{aligned}
\end{equation}
Since 
$$
C:=\int_{R_0}^{R_2} \frac { \int_0^r\sqrt{\rho_1(s)}\,ds}{\psi(r)+\log(\rho_2(r))}\sqrt{\rho_1(r)}\,dr<+\infty,
$$
then we have that $\varphi$ satisfies \eqref{integral-varphi}. Theorem \ref{teo1} yields the thesis.

\section{Proof of Theorems \ref{main1-p1}, \ref{main2-p1} and  \ref{teogen-p1}}\label{sec5}\setcounter{equation}{0}

Similarly to the previous case, in order to prove Theorems \ref{main1-p1}, \ref{main2-p1} and  \ref{teogen-p1}, we take advantage of a general uniqueness result for the Cauchy problem of the heat equation in \eqref{problemachanged} where $M$ is a weighted, geodetically complete Riemannian manifold. The goal will be, by means of a suitable conformal change of the metric defined in \eqref{conformalchange}, to reformulate problem \eqref{problema} into problem \eqref{problemachanged} where the operator is given in \eqref{deltaweighted} and the measure is \eqref{weighted_mess}.

\subsection{Uniqueness for the heat equation on a weighted manifold with $p=1$}

We first need to introduce the following assumption which will restrict the class of manifolds that we can consider. 
\smallskip

We will say that the weighted manifold $(M,g,\mu_\omega)$ satisfies assumption \eqref{acca} if
\begin{equation}\label{acca}\tag{H}
\begin{aligned}
& \,\text{(i)}\,\,\text{$M$ has a pole $o$;}\\
&\,\text{(ii)}\,\,\Delta r\ge0\,\,\,\text{where}\,\,\,r=\operatorname{dist}(x,o);\\
&\,\text{(iii)}\,\,\text{$\omega\in C^\infty(M)$, $\omega=\omega(r)$, $\omega>0$, $\frac{\partial\omega}{\partial r}= O(r)$} \text{ for every } r>1\,.
\end{aligned}
\end{equation}
The following theorem is formulated for any weighted, geodetically complete Riemannian manifold. In this subsection, the weight $e^{-\omega}$ does not depend on $\rho$ (where we have used the notation introduced in Section \ref{sec3}). For this reason, the following theorem is of independent interest.

\begin{theorem}\label{teo2}
Let $(M, g,\mu_{\omega})$ be a weighted Riemannian manifold such that \eqref{acca} holds. Moreover, let $\varphi:M\to(0,+\infty)$ be an increasing and continuous function such that, for some $x_0\in M$ and $R_0>0$,
\begin{equation}\label{peso}
\int_{R_0}^\infty \frac{ r}{\varphi(r)}\,d r=+\infty\,.
\end{equation}
Let $u$ be a solution to problem \eqref{problemachanged} such that
\begin{equation}\label{mainassumptiontilde2}
\int_0^T\int_{B_R(x_0)}|u(x,t)|\,d\mu(x)dt\le e^{\varphi(R)}\quad\text{for all}\,\,R> R_0\,.
\end{equation}
Then 
$$
u\equiv 0 \quad \text{ in } \quad M\times(0,T]\, .
$$
\end{theorem}

A similar result in the unweighted setting has already been obtained in \cite[Theorem 2.4]{Punzo2015}. Although, due to the weighted setting and since our assumptions are different, we provide the proof of Theorem \ref{teo2} in the next subsection.

\

\noindent \textbf{Proof of Theorem \ref{teo2}.} In the proof of Theorem \ref{teo2}, we will take advantage of the following Lemma which has been proved already in \cite[Theorem 9.2]{Grig} for $p=2$ and generalized to $p=1$ in \cite[Lemma 3.1]{Punzo2015}. 

\begin{lemma}\label{lemma2} 
Let assumptions of Theorem \ref{teo2} be satisfied. Moreover, assume that there exist $C>0$, $R_0>0$, $\tau\in(0,T)$ and $x_0\in M$ such that
$$
\int_{\tau-\delta}^\tau\int_{B_R(x_0)}|u(x,\tau)|\,d\mu_\omega(x)dt\le\int_{\tau-\delta}^\tau \int_{B_{2R}(x_0)}|u(x,\tau-\delta)|\,d\mu_\omega(x)dt\,+\frac{C}{R^2}\,,
$$
for any $R>R_0$ and for any $\delta\in\left(0,\min\left\{\tau,\frac{R^2}{C\varphi(2R)}\right\}\right)$, where $\varphi$ is the same as Theorem \ref{teo2}. Then 
$$
u\equiv 0 \quad \text{ in } \quad M\times(0,T]\, .
$$
\end{lemma}

We now provide an a priori estimate which involves a cut-off function and a regular enough test function. Then the theorem will follow by properly construct such two functions.
\begin{proposition}\label{prop2}
Let $(M, g, \mu_\omega)$ be a complete, noncompact, Riemannian manifold such that \eqref{acca} holds. Let $u$ be a solution to problem \eqref{problemachanged}. Let $\eta\in C_c^\infty(M)$ be a family of cut-off functions. Let $R>0$ and be given a function $\zeta\in C^1(M\times(0,T])$, with $\zeta(\cdot,t)\in C^2(M\setminus\partial B_R)$ for all $t\in(0,T]$. 
Then for $\tau\in(0,T)$ and $\delta\in(0,\tau)$
$$
\begin{aligned}
\int_M |u(x,\tau)|\eta^2(x)\,e^{\zeta(x,\tau)}d\tilde\mu_\omega(x)&\le\int_M |u(x,\tau-\delta)|\eta^2(x)\,e^{\zeta(x,\tau-\delta)}d\mu_\omega(x)\\
&+ 2 \int_{\tau-\delta}^\tau\int_{M}|u(x,t)| \,e^{\zeta}\left(\eta\Delta_{\omega} \eta+2|\nabla\eta|^2\right)d\mu_\omega(x) dt\\
&+ \int_{\tau-\delta}^\tau\int_{M\setminus B_R}|u(x,t)|\,\eta^2e^\zeta\left(\partial_t\zeta+3|\nabla \zeta|^2+\Delta_{\omega} \zeta\right)d \mu_\omega(x) dt\\
&+ \int_{\tau-\delta}^\tau\int_{B_R}|u(x,t)|\eta^2e^\zeta\left(\partial_t\zeta+3|\nabla\zeta|^2+\Delta_\omega \zeta\right)d \mu_\omega(x) dt\, . 
\end{aligned}
$$
\end{proposition}

\begin{proof}
For any $\alpha>0$ small enough, we define
\begin{equation}\label{G}
G_{\alpha}(u):=(u^2+\alpha)^{\frac 12}\,.
\end{equation}
We consider expression 
$$\partial_t [G_\alpha(u)]-\Delta_{\omega}[G_{\alpha}(u)], $$
we multiply it by $\eta^2\,e^{\zeta}$ and we integrate over $M\times(\tau-\delta,\tau)$, thus we write
\begin{equation}\label{eq61}
\int_{\tau-\delta}^\tau\int_M \left\{\partial_t [G_\alpha(u)]\,\eta^2\,e^{\zeta}-\Delta_{\omega}[G_{\alpha}(u)]\eta^2\,e^{\zeta}\right\}\,d\mu_\omega dt\,.
\end{equation}
By exploiting the derivatives in \eqref{eq61} and due to problem \eqref{problemachanged}, we get
\begin{equation}\label{eq61b}
\begin{aligned}
\int_{\tau-\delta}^\tau\int_M& \left\{\partial_t [G_\alpha(u)]\,\eta^2\,e^{\zeta}-\Delta_\omega[G_{\alpha}(u)]\eta^2\,e^{\zeta}\right\}\,d\mu_\omega dt\\
&=\int_{\tau-\delta}^\tau\int_M \left\{ [G_\alpha(u)]'\left[\partial_tu-\Delta_\omega u\right]-[G_{\alpha}(u)]''|\nabla u|\right\}\eta^2\,e^{\zeta}\,d\mu_\omega dt\\
&=-\int_{\tau-\delta}^\tau\int_M [G_{\alpha}(u)]''|\nabla u|\eta^2\,e^{\zeta}\,d\mu_\omega dt\,.\end{aligned}
\end{equation}
Observe that for each $\alpha>0$
$$
(G_\alpha)''(t)\ge0 \quad\text{for all}\,\,\,t\in\R\,,
$$
therefore, from \eqref{eq61} and \eqref{eq61b}, we get
\begin{equation}\label{eq62}
\int_{\tau-\delta}^\tau\int_M \partial_t [G_\alpha(u)]\,\eta^2\,e^{\zeta}\,d\mu_\omega dt\le \int_{\tau-\delta}^\tau\int_M \Delta_\omega[G_{\alpha}(u)]\eta^2\,e^{\zeta}\,d\mu_\omega dt\,.
\end{equation}
If we integrate by parts the left-hand-side of \eqref{eq62}, we obtain
\begin{equation}\label{eq63}
\begin{aligned}
\int_{\tau-\delta}^\tau\int_M &\partial_t [G_\alpha(u)]\,\eta^2\,e^{\zeta}\,d\mu_\omega dt\\
&=\partial_t\left\{\int_{\tau-\delta}^\tau\int_M [G_\alpha(u)]\,\eta^2\,e^{\zeta}\,d\mu_\omega dt\right\}-\int_{\tau-\delta}^\tau\int_M [G_\alpha(u)]\,\eta^2\,\partial_t\left(e^{\zeta}\right)\,d\mu_\omega dt\\
&=\int_M [G_\alpha(u(x,\tau))]\,\eta^2\,e^{\zeta(x,\tau)}\,d\mu_\omega(x) -\int_M [G_\alpha(u(x,\tau-\delta))]\,\eta^2\,e^{\zeta(x,\tau-\delta)}\,d\mu_\omega(x)\\
&\quad -\int_{\tau-\delta}^\tau\int_M [G_\alpha(u)]\,\eta^2e^{\zeta}\,\partial_t\zeta \,d\mu_\omega dt\, .
\end{aligned}
\end{equation}
From \eqref{eq62} and \eqref{eq63} we obtain
\begin{equation}\label{eq64}
\begin{aligned}
\int_M [G_\alpha(u(x,\tau))]&\,\eta^2e^{\zeta(x,\tau)}\,d\mu_\omega\le \int_M [G_\alpha(u(x,\tau-\delta))]\,\eta^2e^{\zeta(x,\tau-\delta)}\,d\mu_\omega(x)\\
& +\int_{\tau-\delta}^\tau\int_M [G_\alpha(u)]\,\eta^2 e^{\zeta}\,\partial_\omega\zeta\,d\mu_f(x) dt + \int_{\tau-\delta}^\tau\int_M \Delta_\omega[G_{\alpha}(u)]\eta^2e^{\zeta}\,d\mu_\omega dt
\end{aligned}
\end{equation}
Since $\zeta(\cdot,t)\in C^2(M\setminus\partial B_R)$ for all $t\in(0,T)$, then for any $t\in(\tau-\delta,\tau)$, we have
\begin{equation}\label{eq65}
\begin{aligned}
\int_M &\Delta_{\omega}[G_{\alpha}(u)]\eta^2e^{\zeta}\,d\mu_\omega =\int_{B_R} \Delta_{\omega}[G_{\alpha}(u)]\eta^2e^{\zeta}\,d\mu_\omega +\int_{M\setminus B_R} \Delta_{\omega}[G_{\alpha}(u)]\eta^2e^{\zeta}\,d\mu_\omega \\
=&-\int_{B_R} g\left(\nabla [G_{\alpha}(u)],\nabla\left(\eta^2e^{\zeta}\right)\right)\,d\mu_\omega-\int_{M\setminus \tilde B_R} g\,\left(\nabla[G_{\alpha}(u)],\nabla \left(\eta^2e^{\zeta}\right)\right)\,d\mu_\omega\\
&+\int_{\partial \tilde B_R}\eta^2e^{\zeta}\partial_n G_{\alpha}(u)\,d S_\omega-\int_{\partial B_R}\eta^2e^{\zeta}\partial_n\left(G_{\alpha}(u)\right)\,d S_\omega\\
=&\int_{B_R}G_{\alpha}(u) \Delta_{\omega} \left(\eta^2e^{\zeta}\right)d\mu_\omega-\int_{\partial B_R}G_{\alpha}(u)\partial_n\left(\eta^2 e^\zeta\right)\,d S_\omega\\
&\,\,+\int_{M\setminus B_R}G_{\alpha}(u) \Delta_{\omega} \left(\eta^2e^{\zeta}\right)d\mu_\omega +\int_{\partial B_R}G_{\alpha}(u)\partial_n\left(\eta^2 e^\zeta\right)\,d S_\omega\\
=&2\int_{M}G_{\alpha}(u) e^{\zeta}\left(\eta \Delta_{\omega} \eta+|\nabla\eta|^2\right)\,d\mu_\omega+4\int_{M}G_{\alpha}(u)\eta e^{\zeta}\, g\left(\nabla\eta,\,\nabla \zeta\right)\,d\mu_\omega\\
&\quad+\int_{M\setminus B_R}G_{\alpha}(u)\eta^2e^\zeta\left(|\nabla \zeta|^2+\Delta_{\omega} \zeta\right)d \mu_\omega+\int_{ B_R}G_{\alpha}(u)\eta^2e^\zeta\left(|\nabla\zeta|^2+\Delta_{\omega} \zeta\right)d \mu_\omega\,. 
\end{aligned}
\end{equation}
By means of Cauhy-Schwarz and Young inequalities,
$$
\begin{aligned}
4\int_{M}G_{\alpha}(u)\,\eta \,e^{\zeta}\, g\left(\nabla\eta,\nabla \zeta\right)\,d\mu_\omega&\le 4\int_{M}G_{\alpha}(u)\,\eta \,e^{\zeta}\,|\nabla\eta||\nabla \zeta|\,d\mu_\omega\\
&\le 2\int_{M}G_{\alpha}(u)\, e^{\zeta}\,|\nabla\eta|^2\,d\mu_\omega+2\int_{M}G_{\alpha}(u)\,\eta^2\, e^{\zeta}|\nabla\zeta|^2\,d\mu_\omega\, .
\end{aligned}
$$
Therefore, \eqref{eq65} reads
\begin{equation}\label{eq66}
\begin{aligned}
\int_M &\Delta_{\omega}[G_{\alpha}(u)]\eta^2\,e^{\zeta}\,d\mu_\omega
\le2\int_{M}G_{\alpha}^2(u) e^{\zeta}\left(\eta \Delta_{\omega} \eta+2|\nabla\eta|^2\right)\,d\mu_\omega\\
&\quad+\int_{M\setminus B_R}G_{\alpha}(u)\eta^2e^\zeta\left(3|\nabla \zeta|^2+\Delta_{\omega} \zeta\right)d \mu_\omega\\
&\quad +\int_{ B_R}G_{\alpha}(u)\eta^2e^\zeta\left(3|\nabla\zeta|^2+\Delta_{\omega} \zeta\right)d \mu_\omega\, .
\end{aligned}
\end{equation}
Due to \eqref{eq66} and \eqref{eq63} we obtain
$$
\begin{aligned}
\int_M [G_\alpha(u(x,\tau))]&\,\eta^2\,e^{\zeta(x,\tau)}\,d\mu_\omega\le \int_M [G_\alpha(u(x,\tau-\delta))]\,\eta^2\,e^{\zeta(x,\tau-\delta)}\,d\mu_\omega\\
&+ 2 \int_{\tau-\delta}^\tau\int_{M}G_{\alpha}(u) e^{\zeta}\left(\eta\Delta_{\omega} \eta+2|\nabla\eta|^2\right)\,d\mu_\omega dt\\
&+ \int_{\tau-\delta}^\tau\int_{M\setminus B_R}G_{\alpha}(u)\,\eta^2e^\zeta\left(\partial_t\zeta+3|\nabla\zeta|^2+\Delta_{\omega} \zeta\right)d \mu_\omega dt\\
&+ \int_{\tau-\delta}^\tau\int_{B_R}G_{\alpha}(u)\eta^2e^\zeta\left(\partial_t\zeta+3|\nabla\zeta|^2+\Delta_{\omega} \zeta\right)d \mu_\omega dt\, . 
\end{aligned}
$$
Sending $\alpha\to0^+$ in the latter, due to the dominated convergence theorem, we get the thesis because
$$
G_\alpha(u)\to |u|\quad\text{as}\,\,\,\alpha\to0^+.
$$
\end{proof}

\begin{lemma}\label{lemma3}
Let $(M, g, \mu_\omega)$ be a complete, noncompact, Riemannian manifold such that \eqref{acca} holds. Let $\delta\in\left(0,\min\{\tau,\frac14\}\right)$.
For every $x\in M$ and $t\in(\tau-\delta,\tau]$, we define, for $\alpha\ge\frac{\lambda}{1-4\delta}$ the function
\begin{equation}\label{xi}
\xi(x,t):=-\frac{[r-R]_+^2}{2\alpha(\tau+\delta-t)}\,.
\end{equation}
Then, for every $x\in M\setminus \partial B_R$ and for every $t\in(\tau-\delta,\tau]$,
\begin{equation}\label{supersolution_pequal1}
\partial_t\xi+\frac{\lambda}{2} |\nabla\xi|^2+\Delta_{\omega}\xi\le0\,.
\end{equation}
\end{lemma}
\begin{proof}
Firstly, we observe that $\xi(t)\in C^2(M\setminus\partial B_R)\cap C^1(M)$ for every $t\in(\tau-\delta,\tau]$. Secondly, if $x\in \overline{B_R}$, then the inequality is trivially verified since $\xi\equiv0$. Finally, let $x\in M\setminus \overline{B_R}$, then by the definition of $\xi $ in \eqref{xi} we have, for every $t\in(\tau-\delta,\tau]$,
$$
\begin{aligned}
&\partial_t \xi(x,t)=-\frac{[r-R]^2}{2\alpha(\tau+\delta-t)^2}\\
&\nabla \xi (x,t)=-\frac{[r-R]}{\alpha(\tau+\delta-t)}\nabla r\quad\text{and}\quad |\nabla \xi |^2(x,t)=\frac{[r-R]^2}{\alpha^2(\tau+\delta-t)^2}|\nabla r|^2.
\end{aligned}
$$
Furthermore, since $\xi $ is a radial function and due to \eqref{lapl_polar} 
$$
\begin{aligned}
\Delta_\omega \xi(x,t)
&=-\frac1{\alpha(\tau+\delta-t)}-\left[m(r,\theta)-\frac{\partial\omega}{\partial r}\right]\frac{r-R}{\alpha(\tau+\delta-t)}.\\
\end{aligned}
$$
Therefore for any $x\in M\setminus \overline{B_R}$ and $t\in(\tau-\delta,\tau]$, due to \eqref{m}, \eqref{acca}-(ii) and \eqref{acca}-(iii)
$$
\begin{aligned}
\partial_t \xi&+\frac{\lambda}{2}  |\nabla \xi |^2+\Delta_{\omega} \xi\\
&=\frac{1}{2\alpha^2(\tau+\delta-t)^2}\left\{-(\alpha-\lambda)[r-R]^2-2\alpha(\tau+\delta-t)-2 \alpha\, \left[m(r,\theta)-\frac{\partial\omega}{\partial r}\right](r-R)(\tau+\delta-t)\right\}\\
&\le\frac{1}{2\alpha^2(\tau+\delta-t)^2}\left\{-(\alpha-\lambda)[r-R]^2-2\alpha(\tau+\delta-t)+2 \alpha\frac{\partial\omega}{\partial r}(r-R)(\tau+\delta-t)\right\}\\
&\le\frac{1}{2\alpha^2(\tau+\delta-t)^2}\left\{-(\alpha-\lambda)[r-R]^2-2\alpha(\tau+\delta-t)+4 \alpha\delta\left|\frac{\partial\omega}{\partial r}\right|(r-R)\right\}\\
&\le\frac{1}{2\alpha^2(\tau+\delta-t)^2}\left\{-(\alpha-\lambda-4\alpha\delta)[r-R]^2-2\alpha(\tau+\delta-t)\right\} \le 0\,,
\end{aligned}
$$
provided $\alpha\ge\frac{\lambda}{1-4\delta}$.
\end{proof}

We are now in position to prove Theorem \ref{teo2}.

\begin{proof}[Proof of Theorem \ref{teo2}]
Let $R>0$, we choose the cut-off function $\tilde\eta\in C^\infty([0,+\infty))$ such that $0\le\tilde\eta\le1$,
\begin{equation}\label{eq71}
|\tilde\eta'|\le1,\quad|\tilde\eta''\\|\le1,\quad\tilde\eta'\le 0\in [0,+\infty),
\end{equation}
and
$$
\tilde\eta(s)=\begin{cases} 1&\quad\text{if}\,\,\,0\le s<\frac32\\
0&\quad \text{if}\,\,\,s>2\,.\end{cases}
$$
Define 
\begin{equation}\label{eta}
\eta_R:=\tilde\eta\left(\frac{r(x)}{R}\right)\quad\text{for all}\,\,\,x\in M\,.
\end{equation}
From \eqref{eq71}, \eqref{m} and \eqref{acca}-(ii) we get, for any $R>1$,
\begin{equation}\label{etagrad}
\begin{aligned}
|\nabla\eta_R|&=\frac1R\left|\tilde\eta'\left(\frac{r(x)}{R}\right)\right|\le\frac1R\quad\text{in}\,\,\,M;   \\
\Delta_\omega\eta_R&=\frac1{R^2}\tilde\eta''\left(\frac{r(x)}{R}\right)+\frac1R\left[m(r,\theta)-\frac{\partial\omega}{\partial r}\right]\tilde\eta'\left(\frac{r(x)}{R}\right)\le \frac1{R^2} \quad\text{in}\,\,\,M.\\
\end{aligned}
\end{equation}
Let $\tau\in(0,T)$; by Proposition \ref{prop2} with $\eta=\eta_R$ and $\zeta=\xi$ we get, for all $\tau\in(0,T)$ and $\delta\in\left(0,\min\{\tau,\frac14\}\right)$ 
$$
\begin{aligned}
\int_M |u(x,\tau)|\eta_R^2(x)\,e^{\xi(x,\tau)}d\mu_\omega(x)\le&\int_M |u(x,\tau-\delta)|\eta_R^2(x)\,e^{\xi(x,\tau-\delta)}d\mu_\omega(x)\\
&+ 2 \int_{\tau-\delta}^\tau\int_{M}|u(x,t)| \,e^{\xi}\left(\eta_R\Delta_{\omega} \eta_R+2|\nabla\eta_R|^2\right)d\mu_\omega(x) dt\\
&+ \int_{\tau-\delta}^\tau\int_{M\setminus B_R}|u(x,t)|\,\eta_R^2e^ \xi\left(\partial_t \xi +3|\nabla \xi |^2+\Delta_{\omega} \xi\right)d \mu_\omega(x) dt\\
&+ \int_{\tau-\delta}^\tau\int_{ \tilde B_R}|u(x,t)|\eta_R^2e^ \xi\left(\partial_t \xi +3| \nabla\xi |^2+ \Delta_{\omega} \xi\right)d \mu_\omega(x) dt
\end{aligned}
$$
We can now apply Lemma \ref{lemma3} with $\lambda=6$, therefore
$$
\begin{aligned}
\int_M |u(x,\tau)|\eta_R^2(x)\,e^{\xi(x,\tau)}d\mu_\omega(x)\le&\int_M |u(x,\tau-\delta)|\eta_R^2(x)\,e^{\xi(x,\tau-\delta)}d\mu_\omega(x)\\
&+ 2 \int_{\tau-\delta}^\tau\int_{M}|u(x,t)| \,e^{\xi}\left(\eta_R \Delta_{\omega} \eta_R+2|\nabla\eta_R|^2\right)d\mu_\omega(x) dt.
\end{aligned}
$$
Furthermore, due to \eqref{etagrad} and the definition of $\xi$ in \eqref{xi}, the latter reduces to
\begin{equation}\label{eq72}
\begin{aligned}
\int_{\tilde B_{R}}|u(x,\tau)|\,d\mu_\omega(x)&\le\int_{B_{2R}} |u(x,\tau-\delta)|\,d\mu_\omega(x)\\
&+ \frac{2}{R^2} \int_{\tau-\delta}^\tau\int_{B_{2R}\setminus B_{\frac32 R}}|u(x,t)| \,e^{\xi}d\mu_\omega(x) dt.
\end{aligned}
\end{equation}
Moreover, observe that, for any $x\in B_{2R}(x_0)\setminus B_{\frac32R}(x_0)$ and $t\in(\tau-\delta, \tau)$, we have
$$
e^{\xi(x,t)}=\operatorname{exp}\left\{-\frac{[r(x)-R]_+^2}{2\alpha(\tau+\delta-t)}\right\}\le \operatorname{exp}\left\{-\frac{R^2}{16\alpha\delta}\right\}\,.
$$
By plugging it into \eqref{eq72}, we obtain
$$
\begin{aligned}
\int_{B_{R}}|u(x,\tau)|\,d\mu_\omega(x)&\le\int_{B_{2R}} |u(x,\tau-\delta)|\,d\mu_\omega(x)+ \frac{2}{R^2} \int_{\tau-\delta}^\tau\int_{B_{2R}\setminus B_{\frac32 R}}|u(x,t)| \,e^{-\frac{R^2}{16\alpha\delta}}d\mu_\omega(x) dt.
\end{aligned}
$$
Finally, by assumption \eqref{mainassumptiontilde2} of Theorem \ref{teo2}, we write
$$
\begin{aligned}
\int_{B_{R}(x_0)} |u(x,\tau)|d\mu_\omega(x)&\le\int_{B_{2R}(x_0)} |u(x,\tau-\delta)|\,d\mu_\omega(x)+ \frac{2}{R^2}e^{\varphi(2R)-\frac{R^2}{16\alpha\delta}}\,.
\end{aligned}
$$
For any $R> R_0$ and $\delta\in\left(0,\min\left\{\tau,\frac{R^2}{C\varphi(2R)},\frac14\right\}\right)$, there holds, for some $\tilde C>0$
$$
\begin{aligned}
\int_{B_{R}(x_0)} |u(x,\tau)|d\mu_\omega(x)&\le\int_{B_{2R}(x_0)} |u(x,\tau-\delta)|\,d\mu_\omega(x)+ \frac{\tilde C}{R^2}\,.
\end{aligned}
$$
Therefore the thesis follows by means of Lemma \ref{lemma2}.
\end{proof}

\bigskip

\noindent The proofs of Theorems  \ref{main1-p1}, \ref{main2-p1} and  \ref{teogen-p1} are based on the same conformal change of the metric made in Section \ref{sec4}, see \eqref{eq75}. Therefore, we work with the Cauchy problem of the heat equation defined in \eqref{problema:support} and we show that Theorem \ref{teo2} applies. 

\

\noindent \textbf{Proof of Theorem \ref{main1-p1}.} 
Firstly, observe that
$$
\tilde\Delta \tilde r=\frac1{\sqrt{\rho(r)}}\left(\Delta r+\frac{N-1}{2}\frac{\rho'(r)}{\rho(r)}\right)\,.
$$
Then, by means of \eqref{hp-laplaciano-theta}, \eqref{acca}-(ii) follows. Recall that, by \eqref{eq75}, $\omega=\log\frac{1}{\rho}$ therefore 
$$
\frac{\partial\omega}{\partial r}=\frac{\theta r}{1+r}=O(r)\quad \text{for all}\,\,r>1\,.
$$
This ensures that \eqref{acca}-(iii) holds true.
Finally, arguing as in the proof of Theorem \ref{main1} with $p=1$ we get that assumptions \eqref{peso} and \eqref{mainassumptiontilde2} of Theorem \ref{teo2} hold with $\varphi$ defined as in \eqref{eq:phi}.  The thesis follows by Theorem \ref{teo2}.

\

\noindent \textbf{Proof of Theorem \ref{main2-p1}.} 
Firstly, observe that
$$
\tilde\Delta \tilde r=\frac1{\sqrt{\rho(r)}}\left(\Delta r+\frac{N-1}{2}\frac{\rho'(r)}{\rho(r)}\right)\,.
$$
Then, by means of \eqref{hp-laplaciano}, \eqref{acca}-(ii) follows.  Recall that, by \eqref{eq75}, $\omega=\log\frac{1}{\rho}$ therefore 
$$
\frac{\partial\omega}{\partial r}=\frac{2 r}{1+r}=O(r)\quad \text{for all}\,\,r>1\,.
$$
This ensures that \eqref{acca}-(iii) holds true.
Finally, arguing as in the proof of Theorem \ref{main2} with $p=1$ we get that assumptions \eqref{peso} and \eqref{mainassumptiontilde2} of Theorem \ref{teo2} hold with $\varphi$ defined as in \eqref{eq:phi} with $\theta=2$.  The thesis follows by Theorem \ref{teo2}.

\

\noindent \textbf{Proof of Theorem \ref{teogen-p1}.}  
Firstly, observe that
$$
\tilde\Delta \tilde r=\frac1{\sqrt{\rho(r)}}\left(\Delta r+\frac{N-1}{2}\frac{\rho'(r)}{\rho(r)}\right)\,.
$$
Then, by means of \eqref{hp-laplaciano-p1}, \eqref{acca}-(ii) follows. Recall that, by \eqref{eq75}, $\omega=\log\frac{1}{\rho}$ therefore, by \eqref{decayrho}
$$
\frac{\partial\omega}{\partial r}=-\frac{\rho'(r)}{\rho(r)}=O(r)\quad \text{for all}\,\,r>1\,.
$$
It follows that \eqref{acca}-(iii) holds true.
Finally, arguing as in the proof of Theorem \ref{teogen} with $p=1$ we get that assumptions \eqref{peso} and \eqref{mainassumptiontilde2} of Theorem \ref{teo2} hold with $\varphi$ defined as in \eqref{eq:phi-gen}.  The thesis follows by Theorem \ref{teo2}.

\section{General non-uniqueness criteria on model manifolds}\setcounter{equation}{0}\label{sec-nonuniq}

We now show a general non-uniqueness criteria for problem \eqref{problema} in the special case of $M$ being a model manifold, i.e. $M=\mathbb{M}^N_f$ (see Definition \ref{def_modello}). We assume through this section that \eqref{eq:rhobound-general} holds.

\begin{proposition}\label{prop-nonuniq}
Let $\mathbb{M}^N_f$ be a model manifold and let $\rho$ be as in \eqref{eq:rhobound-general}. Suppose that there exist a supersolution $h>0$ of equation 
\begin{equation}\label{eq100}
\Delta h=-\rho(x) \quad\quad \text{in}\,\,\, \mathbb{M}^N_f\,,
\end{equation}
such that 
\begin{equation}\label{eq542}
 \lim_{r(x)\to +\infty} h(x)=0\,.
\end{equation}
Then there exist infinitely many bounded solutions $u$ of problem \eqref{problema0} with $f\equiv0$ and $u_0\in L^\infty(\mathbb M_f^N)$. In particular, for any $\gamma>0$ there exists a solution $u$ to equation \eqref{problema} such that
$$
\lim_{r(x)\to\infty}\frac1T\int_0^T u(x,t)\,dt=\gamma.
$$
\end{proposition}

%

To prove Proposition \ref{prop-nonuniq}, we introduce the following definitions. Let $l=\{1,-1\}$, for some $\gamma>0$, $\chi>0$, we consider the following auxiliary problems:
\begin{equation}\label{CAP2303}
\begin{cases}
\rho u_t=\Delta u\quad &\text{in}\,\, B_R\times(0,T]\\
u=\gamma &\text{on}\,\, \partial B_R\times(0,T]\\
u=u_0 \,\,\,\,\,&\text{on}\,\, B_R\times\{0\},
\end{cases}
\end{equation}
and
\begin{equation}\label{CAP2304}
\begin{cases}
 \Delta V=\,l\,\rho \quad& \text{in}\,\, B_R\\
V=\chi \quad & \text{on}\,\, \partial B_R.
\end{cases}
\end{equation}
Solutions to problems \eqref{CAP2303} and \eqref{CAP2304} are defined as follows.

\begin{definition}\label{CAP2305}
By a solution to the problem \eqref{CAP2303} we mean any function $u\in C(B_R\times(0,T])$ such that
\begin{equation*}
\begin{aligned}
&\int_0^t\int_{B_R}\{\rho u\psi_t+u\Delta \psi\}\,d\mu d\tau=\int_{B_R} \rho\{u(t)\psi(t)-u_0\psi(0)\}\,d\mu+\int_{0}^{t}\int_{\partial B_R}\gamma\,g(\nabla\psi,\nu)\,d\sigma d\tau
\end{aligned}
\end{equation*}
for any $t\in(0,T]$ and any $\psi\in C^{\infty}(\overline{B_R}\times(0,T])$, $\psi\ge 0$ and $\psi =0$ on $\partial B_R\times(0,T]$. Subsolutions (supersolutions) of \eqref{CAP2303} are defined replacing $"="$ by $"\ge"$ (respectively $"\le"$).
\end{definition}

\begin{definition}\label{CAP2307}
By a solution to the problem \eqref{CAP2304} we mean any function $V\in C(\overline{B_R})$ such that
\begin{equation*}
\int_{B_R}V\Delta \psi\,d\mu=\int_{\partial B_R}\chi\, g(\nabla\psi,\nu)\,d\sigma+l\int_{B_R} \rho(x)\psi\,d\mu
\end{equation*}
for any $\psi\in C^{\infty}(\overline{B_R})$, $\psi\ge 0$ and $\psi=0$ on $\partial B_R$. Subsolutions (supersolutions) of \eqref{CAP2304} are defined replacing $"="$ by $"\ge"$ (respectively $"\le"$).
\end{definition}

\begin{proof}[Proof of Proposition \ref{prop-nonuniq}]
 By classical results, we consider a sequence of nonnegative bounded solutions $\{u_{R}\}$ to problem \eqref{CAP2303}. Observe that, since $u_0$ is bounded in $\mathbb M_f^N$, there exists a $K$ independent of $R$ such that for all $(x,t)\in B_R\times(0,T)$
$$
u_R(x,t)\le K\quad \text{for all}\,\,\,R>0.
$$ 
Moreover
thanks to compactness argument, we can extract a subsequence $\{u_{R_j}\}$ where $R_j\to +\infty$ as $j\to +\infty$ that converges  to $u$, which is a solution to problem \eqref{problema} in the sense of Definition \ref{defsole}. Finally, it is still true that
\begin{equation}\label{CAP2309}
0\le u\le K \quad\,\, \text{in}\,\,\mathbb{M}^N_f\times(0,T].
\end{equation}
It remains to show that 
$$
\lim_{r(x)\to\infty} \frac{1}{T}\int_0^Tu(x,t)\,dt=\gamma.
$$
For simplicity of notation, we denote the subsequence again as $u_R$. Define,
\begin{equation}\label{CAP2313}
v_{R}(x):=\int_0^T u_{R}(x,t)\,dt\, , \text{ for any } x\in B_R.
\end{equation}
For any $\psi\in C^\infty(\overline{B_R})$, $\psi=0$ on $\partial B_R$, we easily obtain from Definition \ref{CAP2305}
\begin{equation}\label{2314eq}
\int_{B_{R}}v_{R}\Delta\psi\,d\mu=\int_{B_{R}}\rho\left[u_{R}(x,T)-u_0(x)\right]\psi\,d\mu+ \int_{\partial B_{R}}\gamma\,T\,g(\nabla\psi,\nu)\,d\sigma\, . 
\end{equation}
Moreover, observe that, thanks to \eqref{CAP2309}
\begin{equation}\label{CAP23014}
\left|u_{R}(x,T)\right|+\left|u_0(x)\right|\le 2\,K=:M\,.
\end{equation}
Thus,
\begin{equation}\label{CAP2315}
\begin{aligned}
\int_{B_{R}}v_{R}\Delta\psi\,d\mu&\ge\int_{\partial B_{R}}\gamma T\, g(\nabla\psi,\nu)\,d\sigma -M\int_{B_{R}}\rho\,\psi\,d\mu
\end{aligned}
\end{equation}
Inequality \eqref{CAP2315} shows, dividing by $M$, that the function
\begin{equation}\label{CAP2316}
F_{1}:=\frac{v_{R}}{M},
\end{equation}
is a subsolution of problem \eqref{CAP2304} with 
$l=-1$ and  $\chi:=\dfrac{\gamma\,T}{M}.$
Furthermore, by assumption, there exists $h$ such that
\begin{equation}\label{CAP2318}
\int_{B_{R}}h\Delta\psi\,d\mu\,\le -\int_{B_{R}}\rho\psi\,d\mu+\int_{\partial B_{R}}h\, g(\nabla\psi,\nu)\,d\sigma
\end{equation}
where we have used that $h>0$. We now consider the constant solution $W=\gamma$ of the problem 
$$
\begin{cases}
\Delta W=0 \quad &\text{in}\,\,B_R\\
W=\gamma \quad\quad &\text{on}\,\,\partial B_R.
\end{cases}
$$
Then, for any $\psi\in C^{\infty}(\overline{B_R})$, $\psi\ge 0$ and $\psi=0$ on $\partial B_R$, it follows that
\begin{equation}\label{CAP2319}
\int_{B_R}W\Delta \psi\,d\mu=\int_{\partial B_R}W\, g(\nabla\psi,\nu)\,d\sigma.
\end{equation}
We multiply \eqref{CAP2319} by $\dfrac{T}{M}$ and we sum the result together with \eqref{CAP2318}. Using the definition of $W$ we get,
\begin{equation}\label{CAP2320}
\int_{B_{R}}\left(h+\frac{\gamma\,T}{M}\right)\Delta\psi\,d\mu\,\le\int_{\partial B_{R}}\left(h+\frac{\gamma\,T}{M}\right)g(\nabla\psi,\nu)\,d\sigma -\int_{B_{R}}\rho\psi\,d\mu.
\end{equation}
Defining,
\begin{equation}\label{CAP2317}
F_{2}:=h+\frac{\gamma\,T}{M},
\end{equation}
inequality \eqref{CAP2320} becomes,
\begin{equation}\label{CAP2321}
\int_{B_{R}}F_{2}\Delta\psi\,d\mu\,\le\int_{\partial B_{R}}F_{2}\, g(\nabla\psi,\nu)\,d\sigma-\int_{B_{R}}\rho\psi\,d\mu.
\end{equation}
This proves that $F_{2}$ is a supersolution to problem \eqref{CAP2304}. By comparison results, it follows that
$$
F_{1}\le F_{2}.
$$
Hence
$$
v_{R}=M\,F_{1}\le M\,F_{2}=M\left(h+\frac{\gamma\,T}{M}\right)=M\,h+\gamma\,T \quad \text{in}\,\,B_{R}.
$$
Letting $R\longrightarrow \infty$ we obtain,
\begin{equation}\label{CAP2322}
v\le M\,h+\gamma\,T \quad \text{in}\,\,\mathbb M_f^N,
\end{equation}
where
\begin{equation}\label{CAP2322b}
\begin{aligned}
u&:=\lim_{R\to\infty}u_{R} \quad \text{in} \,\,\mathbb M_f^N\times(0,T],\\
v&:=\int_{0}^{T}u(x,t)\,dt=\int_0^T\lim_{R\to+\infty} u_R(x,t)\,dt=\lim_{R\to+\infty}v_{R} \quad \text{for}\,\, x\in \mathbb M_f^N.
\end{aligned}
\end{equation}
On the other hand, thanks to \eqref{CAP2309}, \eqref{2314eq} and \eqref{CAP23014}
\begin{equation}\label{CAP2323}
\begin{aligned}
\int_{B_{R}}v_{R}\Delta\psi\,d\mu&\le\int_{\partial B_{R}}\gamma T\,g(\nabla\psi,\nu)\,d\sigma+M\int_{B_{R}}\rho\,\psi\,d\mu
\end{aligned}
\end{equation}
Inequality \eqref{CAP2323} shows that the function $F_{1}$ defined in \eqref{CAP2316}, is a supersolution of problem \eqref{CAP2304} for 
$l=1$ and $\chi:=\frac{\gamma\,T}{M}.$
Again, by assumption, we have
\begin{equation}\label{CAP2324}
\int_{B_{R}}(-h)\Delta\psi\,d\mu\,\ge \int_{B_{R}}\rho\psi\,dx+\int_{\partial B_{R}}(-h)g(\nabla\psi,\nu)\,d\sigma,
\end{equation}
where we have used that $h\ge0$.
We now consider the constant solution $W=\gamma$ as defined in \eqref{CAP2319}. Then, we multiply \eqref{CAP2319} by $\dfrac{T}{M}$ and we sum the result together with \eqref{CAP2324}. Using the definition of $W$ we get,
\begin{equation}\label{CAP2325}
\int_{B_{R}}\left(-h+\frac{\gamma\,T}{M}\right)\Delta\psi\,d\mu\,\ge\int_{\partial B_{R}}\left(-h+\frac{\gamma\,T}{M}\right)g(\nabla\psi,\nu)\,d\sigma+\int_{\Omega}\rho\psi\,d\mu.
\end{equation}
Defining
\begin{equation}\label{326}
F_{3}:=-h+\frac{\gamma\,T}{M},
\end{equation}
inequality \eqref{CAP2325} proves that $F_{3}$ is a subsolution to problem \eqref{CAP2304} with 
$l=1$ and  $\chi:=\frac{\gamma\,T}{M}.$
By comparison results, it follows that
$$
F_{1}\ge F_{3},
$$
hence
$$
v_{R}=M\,F_{1}\ge M\,F_{3}=M\left[-h+\frac{\gamma\,T}{M}\right]=-M\,h+\gamma\,T \quad \text{in}\,\,B_{R}.
$$
Letting $R\longrightarrow \infty$ we obtain,
\begin{equation}\label{CAP2328}
v\ge -M\,h+\gamma\,T \quad \text{in}\,\,\mathbb M^N_f.
\end{equation}
Combining \eqref{CAP2322} and \eqref{CAP2328}, thanks to the property of $V$, we obtain
\begin{equation}\label{CAP2329}
\gamma\,T\,\,=\,\lim_{r(x)\to+\infty}(-M\,h+\gamma\,T)\,\,\le\lim_{r(x)\to+\infty}v(x)\le\lim_{r(x)\to+\infty}\,(M\,h+\gamma\,T)=\gamma\,T.
\end{equation}
Thus
\begin{equation}\label{CAP2330}
\lim_{r(x)\to+\infty}v(x)\,\,=\gamma\,T,
\end{equation}
Recalling the definition of $v$ and $u$ in \eqref{CAP2322b}, the thesis follows.
\end{proof}

\section{Sharpness of the uniqueness result on model manifolds}\label{sec-examples} \setcounter{equation}{0}

\

\noindent \textbf{Example of uniqueness.} Let us consider problem \eqref{problema} with $M=\mathbb M_f^N$. We prove the following corollaries which vary depending on the volume growth condition on $\mathbb M_f^N$.

\begin{corollary}\label{cor-ex1}
Let $\mathbb M_f^N$ be a model manifold such that
\begin{equation}\label{eq85}
\displaystyle f(r)= e^{\frac{r^\beta}{N-1}},\quad \text{for some}\,\,\,\beta\in(0,2].
\end{equation}
Furthermore, assume that $\rho\in C^\infty(\mathbb M_f^N)$, $\rho>0$, is a radial function given by
\begin{equation}\label{eq86}
\rho(x)\equiv\rho(r)=\left(1+r^2\right)^{-\frac\theta2} \, , \quad \theta\in[0,2-\beta)\,.
\end{equation}
Then there exists at most one bounded solution to problem \eqref{problema}.
\end{corollary}

\begin{proof}[Proof of Corollary \ref{cor-ex1}]
Let $u_1,u_2\in L^\infty(\mathbb M_f^N\times(0,T))$ be two solutions to problem \eqref{problema} with $M=\mathbb M_f^N$. We define $w=u_1-u_2$. Observe that $w\in L^\infty(\mathbb M_f^N\times(0,T))$ and that $w$ solves \eqref{problema} with $M=\mathbb M_f^N$. We show that $w=0$ by means of Theorem \ref{main1}. Trivially, \eqref{eq:rhobound-theta} holds because $\theta\in[0,2-\beta)$ and $\beta\in(0,2]$. Moreover, we may consider $\psi:(0,+\infty)\to \R$, positive, increasing and continuous in $(0,+\infty)$ defined as
\begin{equation}\label{eq83}
\psi(r)= r^{2-\theta} \quad \text{for all}\,\,\,r>0\,.
\end{equation}
Then \eqref{integral-phi-theta} reads as
$$
\int_{R_0}^\infty \frac{r^{1-\theta}}{\psi(r)}\,dr=\int_{R_0}^\infty \frac1{r}\,dr=+\infty\,.
$$
Furthermore, for any $p>1$, $w\in L^p_{e^{-\psi}}(\mathbb M_f^N\times(0,T))$, in fact, due to \eqref{eq84}, \eqref{eq85} and \eqref{eq83},
$$
\begin{aligned}
\int_0^T\int_{\mathbb M_f^N}|w|^p e^{-\psi}\,d\mu dt&\le \|w\|_{L^\infty}^p \int_0^T\int_{\mathbb M_f^N} e^{-\psi}\,d\mu dt\\
&= \|w\|_{L^\infty}^p \int_0^T\iint_{\mathbb S^{N-1}\times(0,\infty)} e^{r^{-2+\theta+\beta}}\,dr\,d\theta\,dt,
\end{aligned}
$$
which converges because $\theta\in[0,2-\beta)$. Therefore, by Theorem \ref{main1}, $w\equiv0$ and the claim follows.
\end{proof}

\begin{remark}\label{rem-ex1}
Arguing similarly to Corollary \ref{cor-ex1}, it is possible to show also uniqueness of unbounded solutions to problem \eqref{problema} with $M=\mathbb M_f^N$ and $f$ as in \eqref{eq85}, provided that $u\in L^p_{e^{-\psi}}(\mathbb M_f^N\times(0,T))$ for some $p>1$, with $\psi$ as in \eqref{eq83}. For example
\begin{itemize}
\item if $u$ is such that, for all $(x,t)\in \mathbb M_f^N\times(0,T)$ 
$$
u(x,t)\le C (1+r)^\alpha\quad \text{for all}\,\,r\ge0,\,\,\,\, \alpha>0,
$$
then $u\in L^p_{e^{-\psi}}(\mathbb M_f^N\times(0,T))$ for all $\rho$ satisfying \eqref{eq86}. 
\item if $u$ is such that, for all $(x,t)\in \mathbb M_f^N\times(0,T)$ 
$$
u(x,t)\le C e^{r^\alpha}\quad \text{for all}\,\,r\ge0,\,\,\,\,\alpha\in (0,2-\theta-\beta),
$$
then $u\in L^p_{e^{-\psi}}(\mathbb M_f^N\times(0,T))$ for $\rho$ satisfying \eqref{eq86}. We furthermore observe that, in the particular case of the Hyperbolic space, i.e. $\beta=1$, we get
$$
\theta\in [0,1)\quad\text{and}\quad \alpha\in(0,1-\theta)\,.
$$
\end{itemize}
\end{remark}

\begin{corollary}\label{cor-ex2}
Let $\mathbb M_f^N$ be a model manifold such that
\begin{equation}\label{eq87}
\displaystyle f(r)=  r^{\frac{\beta}{N-1}} ,\quad \text{for some}\,\,\,\beta>0.
\end{equation}
Furthermore, assume that $\rho\in C^\infty(\mathbb M_f^N)$, $\rho>0$, is a radial function given by
\begin{equation}\label{eq88}
\rho(x)\equiv\rho(r)=\left(1+r^2\right)^{-\frac{\theta}2} \, , \quad \theta\in[0,2]\,.
\end{equation}
Then there exists at most one bounded solution to problem \eqref{problema}.
\end{corollary}

\begin{proof}[Proof of Corollary \ref{cor-ex2}]
Let $u_1,u_2\in L^\infty(\mathbb M_f^N\times(0,T))$ be two solutions to problem \eqref{problema} with $M=\mathbb M_f^N$. We define $w=u_1-u_2$. Observe that $w\in L^\infty(\mathbb M_f^N\times(0,T))$ and that $w$ solves \eqref{problema} with $M=\mathbb M_f^N$. We show that $w=0$ by means of Theorems \ref{main1} and \ref{main2}. 

We first consider the case $\theta\in[0,2)$, therefore \eqref{eq:rhobound-theta} holds. Then we may consider $\psi:(0,+\infty)\to \R$ as in \eqref{eq83}. Then clearly, \eqref{integral-phi-theta} is satisfied. Furthermore, for any $p>1$, $w\in L^p_{e^{-\psi}}(\mathbb M_f^N\times(0,T))$, in fact, due to \eqref{eq84}, \eqref{eq87} and \eqref{eq83},
$$
\begin{aligned}
\int_0^T\int_{\mathbb M_f^N}|w|^p e^{-\psi}\,d\mu dt&\le \|w\|_{L^\infty}^p \int_0^T\int_{\mathbb M_f^N} e^{-\psi}\,d\mu dt\\
&= \|w\|_{L^\infty}^p \int_0^T\iint_{\mathbb S^{N-1}\times(0,\infty)} e^{r^{-2+\theta}} r^\beta\,drd\theta dt,
\end{aligned}
$$
which converges because $\theta\in[0,2)$. Therefore, by Theorem \ref{main1}, $w\equiv0$ and the claim follows.

We now consider $\theta=2$, therefore  \eqref{eq:rhobound} holds. Moreover, we may consider $\psi:(0,+\infty)\to \R$, positive, increasing and continuous in $(0,+\infty)$ defined as
\begin{equation}\label{eq89}
\psi(r)= \log^2\left(r^{\sqrt{\gamma}}\right) \quad \text{for all}\,\,\,r>0\,\,\,\gamma>1\,.
\end{equation}
Then \eqref{integral-phi-theta} reads
$$
\int_{R_0}^\infty \frac{\log r}{r\psi(r)}\,dr=\int_{R_0}^\infty \frac1{{\gamma} r \log r}\,dr=+\infty\,.
$$
Furthermore, $w\in L^p_{e^{-\psi}}(\mathbb M_f^N\times(0,T))$, in fact, due to \eqref{eq84}, \eqref{eq87} and \eqref{eq89},
$$
\begin{aligned}
\int_0^T\int_{\mathbb M_f^N}|w|^p e^{-\psi}\,d\mu dt&\le \|w\|_{L^\infty}^p \int_0^T\int_{\mathbb M_f^N} e^{-\psi}\,d\mu dt\\
&= \|w\|_{L^\infty}^p \int_0^T\iint_{\mathbb S^{N-1}\times(0,\infty)} e^{-\log^2(r^{\sqrt{\gamma}})} r^\beta\,drd\theta dt\\
&= \|w\|_{L^\infty}^p \int_0^T\iint_{\mathbb S^{N-1}\times(0,\infty)} r^{-\gamma \log(r)+\beta}\,drd\theta dt\\
\end{aligned}
$$
which converges for any $\beta>0$. Therefore, by Theorem \ref{main2}, $w\equiv0$ and the claim follows.
\end{proof}

\begin{remark}\label{rem-ex2}
Also in this case, arguing as in Corollary \ref{cor-ex2}, it is possible to show also uniqueness of unbounded solutions to problem \eqref{problema} with $M=\mathbb M_f^N$ and $f$ as in \eqref{eq85}, provided that $u\in L^p_{e^{-\psi}}(\mathbb M_f^N\times(0,T))$ for some $p>1$ and a suitable $\psi$. More precisely, let $\rho$ be as in \eqref{eq88} with $\theta\in[0,2)$ then
\begin{itemize}
\item if $u$ is such that, for all $(x,t)\in \mathbb M_f^N\times(0,T)$ 
$$
u(x,t)\le C (1+r)^\alpha\quad \text{for all}\,\,r\ge0,\,\,\,\, \alpha>0,
$$
then $u\in L^p_{e^{-\psi}}(\mathbb M_f^N\times(0,T))$ with $\psi$ defined in \eqref{eq83}; 
\item if $u$ is such that, for all $(x,t)\in \mathbb M_f^N\times(0,T)$ 
$$
u(x,t)\le C e^{r^\alpha}\quad \text{for all}\,\,r\ge0,\,\,\,\,\alpha\in (0,2-\theta),
$$
then $u\in L^p_{e^{-\psi}}(\mathbb M_f^N\times(0,T))$ with $\psi$ defined in \eqref{eq83}.
\end{itemize}
Let now $\rho$ be as in \eqref{eq88} with $\theta=2$ then, if $u$ is such that, for all $(x,t)\in \mathbb M_f^N\times(0,T)$ 
$$
u(x,t)\le C (1+r)^{\alpha\log r}\quad \text{for all}\,\,r\ge0,\,\,\,\, \alpha\in\left(0,\frac\gamma p\right),
$$
then $u\in L^p_{e^{-\psi}}(\mathbb M_f^N\times(0,T))$ with $\psi$ defined in \eqref{eq89}. Here $\gamma$ is introduced in \eqref{eq89}.
\end{remark}

\

\noindent \textbf{Example of non-uniqueness.} We now show that the assumptions on $\rho$ in Corollaries \ref{cor-ex1} and \ref{cor-ex2}, see \eqref{eq86} and \eqref{eq88}, are sharp. To do so, we take advantage of Proposition \ref{prop-nonuniq} which shows that infinitely many solutions to equation \eqref{problema} exist, if a supersolution $h$ to \eqref{eq100} exists. Therefore, in the following, we aim to show some existence results depending on the model manifold that we are considering. This, combined with Proposition \ref{prop-nonuniq}, implies non-uniqueness.

\begin{corollary}\label{corollary7} 
Let $\mathbb M_f^N$ be a model manifold such that
\begin{equation*}
\displaystyle f(r)= e^{\frac{r^\beta}{N-1}},\quad \text{for some}\,\,\,\beta\in(0,2].
\end{equation*}
Furthermore, assume that $\rho\in C^\infty(\mathbb M_f^N)$, $\rho>0$, is a radial function given by
\begin{equation*}
\rho(x)\equiv\rho(r)=\left(1+r^2\right)^{-\frac\theta2} \, , \quad \theta>2-\beta\,.
\end{equation*}
Then there exist infinitely many bounded solutions $u$ of problem \eqref{problema0} with $f\equiv0$ and $u_0\in L^\infty(\mathbb M_f^N)$. In particular, for any $\gamma\in \mathbb R$, $\gamma>0$ there exists a solution $u$ to equation \eqref{problema} such that
$$
\lim_{r(x)\to\infty}\frac1T\int_0^T u(x,t)\,dt=\gamma.
$$
\end{corollary}

To prove this we first need to show the following lemma which ensure the existence of a supersolution to problem \eqref{eq100} such as the one needed in Proposition \ref{prop-nonuniq}.

\begin{lemma}\label{lemma7} 
Let $\mathbb M_f^N$ be a model manifold such that
\begin{equation}\label{eq85b}
\displaystyle f(r)= e^{\frac{r^\beta}{N-1}},\quad \text{for some}\,\,\,\beta\in(0,2].
\end{equation}
Furthermore, assume that $\rho\in C^\infty(\mathbb M_f^N)$, $\rho>0$, is a radial function given by
\begin{equation}\label{eq86b}
\rho(x)\equiv\rho(r)=\left(1+r^2\right)^{-\frac\theta2} \, , \quad \theta>2-\beta\,.
\end{equation}
Then there exists a supersolution $h$ of equation \eqref{eq100} which satisfies \eqref{eq542}.
\end{lemma}

\begin{proof} 
We consider the function
\begin{equation*}\label{eq540}
h(x)\equiv h(r):= \int_r^\infty\frac{\displaystyle \int_0^s\rho(\xi)[f(\xi)+1]^{N-1}\,d\xi}{[f(s)+1]^{N-1}}\,\,ds\quad \text{for any}\,\,x\in \mathbb{M}^N_f.
\end{equation*}
Observe that $h$ is well-defined in $[0,+\infty)$, in fact, due to \eqref{eq85b}, \eqref{eq86b} and by performing the change of variable $\xi=t^{\frac1\beta}$
one has, for some $C=C(\beta)>0$
$$
\begin{aligned}
\displaystyle\int_r^\infty\frac{\displaystyle \int_0^s\rho(\xi)[f(\xi)+1]^{N-1}\,d\xi}{[f(s)+1]^{N-1}}\,\,ds&=\int_r^\infty\frac{\displaystyle \int_0^s\left(1+\xi^2\right)^{-\frac\theta2}\left[e^{\frac{\xi^\beta}{N-1}}+1\right]^{N-1}\,d\xi}{\left[e^{\frac{s^\beta}{N-1}}+1\right]^{N-1}}\,\,ds\\
&=\frac1\beta\int_r^\infty \frac{\displaystyle \int_0^{s^\beta}\left(1+t^{\frac2\beta}\right)^{-\frac\theta2} t^{\frac{2-2\beta}{2\beta}}\left[e^{\frac t{N-1}}+1\right]^{N-1}\,dt}{\left[e^{\frac{s^\beta}{N-1}}+1\right]^{N-1}}\,\,ds\\
&\le C(\beta)\int_r^\infty \frac{\displaystyle \int_0^{s^\beta}(1+t^2)^{\frac{1-\theta-\beta}{2\beta}}e^t\,dt}{e^{s^\beta}}\,\,ds\le C\int_r^\infty s^{1-\theta-\beta}\,ds\,.
\end{aligned}
$$
Then, since $\theta>2-\beta$, we conclude that
$$
\displaystyle\int_r^\infty\frac{\displaystyle \int_0^s\rho(\xi)[f(\xi)+1]^{N-1}\,d\xi}{[f(s)+1]^{N-1}}\,\,ds<+\infty
$$
Moreover, it is easily verified that $h>0$ in $\mathbb{M}^N_f$ and $h\to0$ as $r\to\infty$. Furthermore, we compute
\begin{equation}\label{eq541}
\begin{aligned}
h'(r)&=-\int_0^r\rho(\xi)[f(\xi)+1]^{N-1}\,d\xi\,;\\
h''(r)&=-\rho(r)+(N-1)\frac{f'(r)}{(f(r)+1)^N}\int_0^r\rho(\xi)(f(\xi)+1)^{N-1}\,d\xi\,.
\end{aligned}
\end{equation}
By means of \eqref{Laplaciano_modelli} and due to \eqref{eq541},
$$
\begin{aligned}
\Delta h(x)&=h''(r)+(N-1)\frac{f'(r)}{f(r)}h'(r)\\
&=-\rho(r)+(N-1)\frac{f'(r)}{[f(\xi)+1]^N}\int_0^r\rho(\xi)[f(\xi)+1]^{N-1}\,d\xi\\
&\quad-(N-1)\frac{f'(r)}{f(r)[f(\xi)+1]^{N-1}}\int_0^r\rho(\xi)[f(\xi)+1]^{N-1}\,d\xi\\
&\le-\rho(r)+(N-1)\frac{f'(r)}{[f(\xi)+1]^N}\int_0^r\rho(\xi)[f(\xi)+1]^{N-1}\,d\xi\\
&\quad-(n-1)\frac{f'(r)}{[f(\xi)+1]^{N}}\int_0^r\rho(\xi)[f(\xi)+1]^{N-1}\,d\xi\\
&=-\rho(x)\, , 
\end{aligned}
$$
for all $x\in\mathbb M_f^N$. Hence, we get the thesis.
\end{proof}

\begin{proof}[Proof of Corollary \ref{corollary7}]
Combining Proposition \ref{prop-nonuniq} and Lemma \ref{lemma7} we get the thesis.
\end{proof}

We now consider a model manifold with polynomial volume growth.
\begin{corollary}\label{corollary6} 
Let $\mathbb M_f^N$ be a model manifold such that
\begin{equation*}
\displaystyle f(r)=  r^{\frac{\beta}{N-1}} ,\quad \text{for some}\,\,\,\beta>0.
\end{equation*}
Furthermore, assume that $\rho\in C^\infty(\mathbb M_f^N)$, $\rho>0$, is a radial function given by
\begin{equation*}
\rho(x)\equiv\rho(r)=\left(1+r^2\right)^{-\frac\theta2} \, , \quad \theta>2\,.
\end{equation*}
Then there exist infinitely many bounded solutions $u$ of problem \eqref{problema0} with $f\equiv0$ and $u_0\in L^\infty(\mathbb M_f^N)$. In particular, for any $\gamma>0$ there exists a solution $u$ to equation \eqref{problema} such that
$$
\lim_{r(x)\to\infty}\frac1T\int_0^T u(x,t)\,dt=\gamma.
$$
\end{corollary}

Equivalently to Corollary \ref{corollary7}, to prove the latter we first need to show the following 
\begin{lemma}\label{lemma6} 
Let $\mathbb M_f^N$ be a model manifold such that
\begin{equation}\label{eq87b}
\displaystyle f(r)=  r^{\frac{\beta}{N-1}} ,\quad \text{for some}\,\,\,\beta>0.
\end{equation}
Furthermore, assume that $\rho\in C^\infty(\mathbb M_f^N)$, $\rho>0$, is a radial function given by
\begin{equation}\label{eq88b}
\rho(x)\equiv\rho(r)=\left(1+r^2\right)^{-\frac\theta2} \, , \quad \theta>2\,.
\end{equation}
Then there exists a supersolution $h$ of equation \eqref{eq100} which satisfies \eqref{eq542}.
\end{lemma}

\begin{proof} 
We consider the function
\begin{equation*}\label{eq540}
h(x)\equiv h(r):= \int_r^\infty\frac{\displaystyle \int_0^s\rho(\xi)[f(\xi)+1]^{N-1}\,d\xi}{[f(s)+1]^{N-1}}\,\,ds\quad \text{for any}\,\,x\in \mathbb{M}^N_f.
\end{equation*}
Observe that $h$ is well-defined in $[0,+\infty)$, in fact, arguing as in the proof of Lemma \ref{lemma7}, due to \eqref{eq87b}, \eqref{eq88b} and since $\theta>2$, for some $C=C(\beta)>0$
$$
\int_r^\infty\frac{\displaystyle \int_0^s\rho(\xi)[f(\xi)+1]^{N-1}\,d\xi}{[f(s)+1]^{N-1}}\,\,ds\le C(\beta)\int_r^\infty (1+s^2)^{\frac{-\theta+1}{2}}\,\,ds<+\infty
$$
The rest of the properties follow as in the proof of Lemma \ref{lemma7}.
\end{proof}

\begin{proof}[Proof of Corollary \ref{corollary6}]
Combining Proposition \ref{prop-nonuniq} and Lemma \ref{lemma6} we get the thesis.
\end{proof}

\begin{remark}
Observe that the results in Corollaries \ref{corollary7} and \ref{corollary6} show that our assumptions on the density function $\rho$ made to achieve uniqueness of solutions in Corollaries \ref{cor-ex1} and \ref{cor-ex2}, are sharp. More precisely we have showed that
\begin{itemize}
\item if \begin{equation*}
\displaystyle f(r)= e^{\frac{r^\beta}{N-1}},\quad \text{for some}\,\,\,\beta\in(0,2]\quad\text{and}\quad \rho(r)=\left(1+r^2\right)^{-\frac\theta2}
\end{equation*}
then
\begin{itemize}
\item uniqueness of bounded solutions follows if $\theta\in[0,2-\beta)$;
\item infinitely many bounded solutions exist if $\theta>2-\beta$\,.
\end{itemize}
\item if \begin{equation*}
\displaystyle f(r)=  r^{\frac{\beta}{N-1}} ,\quad \text{for some}\,\,\,\beta>0\quad\text{and}\quad \rho(r)=\left(1+r^2\right)^{-\frac\theta2}
\end{equation*}
then
\begin{itemize}
\item uniqueness of bounded solutions follows if $\theta\in[0,2]$;
\item infinitely many bounded solutions exist if $\theta>2$\,.
\end{itemize}
\end{itemize}
\end{remark}


\

\noindent{\bf Acknowledgement} 

\noindent The first author is funded by the Deutsche Forschungsgemeinschaft (DFG, German Research Foundation) - Project-ID 317210226 - SFB 1283. The second author is funded by the Deutsche Forschungsgemeinschaft (DFG, German Research Foundation) - Project-ID 317210226 - SFB 1283/2 2021. The third author is a member of ``Gruppo Nazionale per le Strutture Algebriche, Geometriche e le loro Applicazioni'' (GNSAGA) of the ``Istituto Nazionale di Alta Matematica'' (INdAM, Italy) and has partially supported by the project PRIN 2022 ``Differential--geometric aspects of manifolds via Global Analysis''.
We are grateful to Fabio Punzo for valuable discussions.
\

\

\noindent{\bf Data availability statement}

\noindent Data sharing not applicable to this article as no datasets were generated or analysed during the current study.

\

\end{document}